\documentclass[twoside, 15pt]{article}
\usepackage{amsmath} 
\usepackage[e]{esvect}
\usepackage{graphicx}
\usepackage{epstopdf}
\usepackage[adobe-utopia]{mathdesign}
\usepackage[T1]{fontenc}
\usepackage[english]{babel}
\usepackage{color}

\usepackage{hyperref}

\title{\textsc{Gluing and cutting cube tiling codes in dimension six}}
\date{}
\author{{Andrzej P. Kisielewicz}\\
{\small \sf  A.Kisielewicz@wmie.uz.zgora.pl}\\
{\small\it Wydzia{\l} Matematyki, Informatyki i Ekonometrii, Uniwersytet Zielonog\'orski,}\\
{\small\it ul. Podg\'orna 50, 65-246 Zielona G\'ora, Poland}}

\newtheorem{lemat}{\sc Lemma}
\newtheorem{tw}{\sc Theorem}
\newtheorem{st}{\sc Statement}
\newtheorem{wn}{\sc Corollary}
\newtheorem{df}{\sc Definition}

\newtheorem{uw}{\sc Remark}
\newtheorem{uwi}[uw]{\sc Remarks}

\newtheorem{nap}{\sc Example}

\newtheorem{nps}[nap]{\sc Examples}

\def\ka #1{\mathscr{#1}}
\def\kal #1 #2{\mathscr{#1}^{#2}}
\def\proof{\noindent \textit{Proof.\,\,\,}}

\def\en{\mathbb{N}}
\def\zet{\mathbb{Z}}

\def\er{\mathbb{R}}

\def\te{\mathbb{T}}

\def\Aut #1 #2{\operatorname{Aut}^{#1} (#2)}

\def\bred #1 {\colorbox{red}{ #1}}
\def\red #1 {{\color{red} #1 }} 

\graphicspath{{figpdfkolor/}}
\DeclareGraphicsExtensions{.pdf,.png} 
\begin{document}
\maketitle


\begin{abstract}
Let $S$ be a set of arbitrary objects, and let $s\mapsto s'$ be a permutation  of $S$ such that $s''=(s')'=s$ and $s'\neq  s$. Let $S^d=\{v_1...v_d\colon v_i\in S\}$.  Two words $v,w\in S^d$ are dichotomous if $v_i=w'_i$ for some $i\in [d]=\{1,...,d\}$, and they form a twin pair if $v_i'=w_i$ and $v_j=w_j$ for every $j\in [d]\setminus \{i\}$. A polybox code is a set $V\subset S^d$ in which every two distinct words are dichotomous. A polybox code $V$ is a cube tiling code if $|V|=2^d$. A $2$-periodic cube tiling of $\er^d$ and a cube tiling of flat torus $\te^d$ can be encoded in a form of a cube tiling code. A twin pair $v,w$ in which $v_i=w_i'$ is glue (at the $i$th position) if the pair $v,w$ is replaced by one word $u$ such that $u_j=v_j=w_j$ for every $j\in [d]\setminus \{i\}$ and $u_i=*$, where $*\not\in S$ is some extra fixed symbol.  A word $u$ with $u_i=*$ is cut (at the $i$th position) if $u$ is replaced by a twin pair $q,t$ such that $q_i=t_i'\neq*$ and $u_j=q_j=t_j$ for every $j\in [d]\setminus \{i\}$. If $V,W\subset S^d$ are two cube tiling codes and there is a sequence of twin pairs which can be interchangeably gluing and cutting in a way which allows us to pass from $V$ to $W$, then we say that $W$ is obtained from $V$ by gluing and cutting. In the paper it is shown that for every two cube tiling codes in dimension six  one can be obtained from the other by gluing and cutting.
\end{abstract}

\section{Introduction}
Let $S$ be a set of arbitrary objects, and let $s\mapsto s'$ be a permutation  of $S$ such that $s''=(s')'=s$ and $s'\neq  s$. In the paper we assume that $S$ is finite. Let $S^d=\{v_1...v_d\colon v_i\in S\}$.  Elements of $S$ are called {\it letters}, while members  of $S^d$ will be called {\it words}. We add to $S$ an extra letter $*$ and the set $S\cup \{*\}$ will be denoted by $*S$. We assume that $*'=*$.  Two words $v,w\in (*S)^d$ are {\it dichotomous} if there is $i\in [d]=\{1,...,d\}$ such that $v_i,w_i\in S$ and $v_i=w'_i$, and if additionally $v_j=w_j$ for $j\in [d]\setminus \{i\}$, then we say that $v$ and $u$ form a {\it twin pair} ({\it in the $i$th direction}).
A {\it polybox code} (or {\it genome}) is a set $V\subset (*S)^d$ in which every two distinct words are dichotomous. (Sometimes we shall write just a `code' instead of `polybox code'.) To examine all non-isomorphic polybox codes $V\subset S^d$ it is enough, by \cite[Lemma 2.2]{Kii}, to assume that $|S|\leq 2^d$, that is  $S=\{a_1,a'_1,...,a_k,a_k'\}$, where $k\leq 2^{d-1}$.  A polybox code $V\subset S^d$ is a {\it cube tiling code} if $|V|=2^d$. It is easy to see that a cube tiling code $V\subset S^d$ induces an $r$-perfect code $C\subset \zet^d_{4r+2}$, $r\in \en$, in the maximum metric (\cite{Co,Kii}).
A natural interpretation of a polybox code is some system of boxes. To describe it we shall give a few definitions. (More information on polybox codes can be found in the first two sections of the papers \cite{Kap,KisL,Kii}.)

Let $X_1,\ldots ,X_d$ be non-empty sets with $|X_i|\geq 2$ for every $i\in [d]$. The set $X=X_1\times\cdots \times X_d$ is called a $d$-{\it box}.
A non-empty set $K \subseteq X$ is called a \textit{ box} if $K=K_1\times\cdots \times K_d$ and
$K_i\subseteq X_i$ for each $i\in [d]$. 
Two boxes $K$ and $G$ in $X$ are called \textit{dichotomous} if there is $i\in [d]$ such that $K_i=X_i\setminus G_i$. A \textit{suit} is any collection of pairwise
dichotomous boxes. 
A non-empty set $F\subseteq X$ is said to be a \textit{ polybox} if
there is a suit $\ka F$ for $F$, that is, if $\bigcup \ka F=F$. In other words, $F$ is a polybox if it has a partition into pairwise dichotomous boxes. 


To pass from polybox codes to suits we shall use a kind of translation of words into boxes: Let $X=X_1\times \cdots \times X_d$ be a $d$-box.
Suppose that for each $i\in[d]$ a mapping $f_i\colon S\to 2^{X_i}\setminus \{\emptyset,X_i\}$ is such that $f_i(s')=X_i\setminus f_i(s)$ (Figure 1  and 2). Additionally, $f_i(*)=X_i$ for $i\in [d]$. We define the mapping $f\colon (*S)^d\to 2^X$ by $f(s_1\ldots s_d)=f_1(s_1)\times\cdots\times f_d(s_d).$ 
If now $V\subset (*S)^d$ is a code, then the set of boxes $f(V)=\{f(v)\colon v\in V\}$ is a suit for the polybox $\bigcup f(V)$. The set $f(V)$ is said to be a \textit{realization} of the code $V$. If $v\in (*S)^d$ and $W\subset (*S)^d$ is a polybox code, then we say that $v$ is {\it covered} by $W$, which is denoted by $v\sqsubseteq W$, if $f(v)\subset \bigcup f(W)$ for every $f$ that preserves dichotomies. If $V\subset (*S)^d$ is a polybox code and $v\sqsubseteq W$ for every $v\in V$, then we say that the code $V$ is {\it covered} by $W$ and write $V\sqsubseteq W$.

A {\it cube tiling} of $\er^d$ is a  family of pairwise disjoint cubes $[0,1)^d+T=\{[0,1)^d+t\colon t\in T\}$ such that $\bigcup_{t\in T}([0,1)^d+t)=\er^d$. 
 A cube tiling $[0,1)^d+T$ is called {\it 2-periodic} if $T+2\zet^d=T$. It is easy to show that any 2-periodic cube tiling of $\er^d$ is a realization of some cube tiling code $V\subset S^d$ (\cite[Section 1]{Kii}). Obviously, a 2-periodic cube tiling $[0,1)^d+T$ of $\er^d$ defines a cube tiling $\ka T=[0,1)^d+T_2$ of the {\it flat torus} $\te^d=\{(x_1,\ldots ,x_d)({\rm mod} 2):(x_1,\ldots ,x_d)\in \er^d\}$, where $T_2=\{(x_1,\ldots ,x_d)({\rm mod} 2):(x_1,\ldots ,x_d)\in T\}$. In general, any realization $f(V)$ of a cube tiling code $V\subset S^d$ is a partition of a $d$-box $X$ into $2^d$ pairwise dichotomous boxes (Figure 1 and 2). Such partition is called a {\it minimal partition} (\cite{GKP,KP}). 

It is easy to check that if $V\subset S^d$ is a polybox code, and $v,w\in V$ is a twin pair such that $v_i=w_i'$, then the set of words $U\subset (*S)^d$ given by $U=V\setminus \{v,w\}\cup \{u\}$, where $u_j=v_j$ for $j\neq i$ and $u_i=*$, is still a polybox code. Indeed, if $U$ is not a polybox code, then there is a word $q \in V\setminus \{v,w\}$ such that $q$ and $u$ are not dichotomous. Thus, (recall that $V$ is a polybox code and $v,w,q\in V$) we have $q_i=v_i'$ and $q_i=w_i'$. But $v_i=w_i'$, and then $q$ and $v$ are not dichotomous words. Therefore, in every polybox code $V$ we can replace a twin pair $v,w\in V$, if $V$ contains such a pair, by the above defined word $u$ obtaining a polybox code $U$. We call the word $u$ a {\it gluing} of $v$ and $w$. Clearly, we may reverse this operation replacing $u$ with $u_i=*$ by a twin pair $q,t$ with $q_j=t_j=u_j$ for $j\neq i$ and $q_i=t'_i\neq *$. This pair $q,t$ will be called a {\it cutting} of $u$. Thus, if $V\subset S^d$ is a polybox code containing twin pairs, then we can obtain a polybox code $W$ from $V$ by gluing a twin pair $v,w\in V$ with $v_i=w_i'$ and next cutting the gluing of $v,w$ obtaining a twin pair $q,t\in W$ with $q_i=t_i'$. We shall say that $W$ is obtained from $V$ by {\it gluing and cutting} if there is a sequence of such local transformations (see Example 1) which lead from $V$ to $W$. More precisely, there are two sequences of twin pairs $(\{v^n,w^n\})_{n=1}^m$ and $(\{q^n,t^n\})_{n=1}^m$ such that $q^k,t^k$ is obtained from $v^k,w^k$ by gluing and next cutting (in the manner described above) for $k\in [m]$, $V^0=V$, $V^k=V^{k-1}\setminus \{v^k,w^k\}\cup \{q^k,t^k\}$, where $\{v^k,w^k\}\subset V^{k-1}$ for $k\in [m]$ and $W=V^m$. (A passing from $V^k$ to $V^{k+1}$ will be called a {\it single operation} in the glue and cut procedure.)

Let us note that gluing and cutting can change a polybox code but not its realization, that is $\bigcup f(V)=\bigcup f(W)$ for every $f$ preserving dichotomies, where $W$ is obtained from $V$ by gluing and cutting. Thus, it is natural to pose the following question: Let $V,W\subset S^d$ be two polybox codes such that $\bigcup f(V)=\bigcup f(W)$ for every $f$ preserving dichotomies (that is $V\sqsubseteq W$ and $W\sqsubseteq V$). We shall call such codes $V,W$ {\it equivalent} (\cite{KP}). 

\medskip
{\bf Question 1}. {\it For which equivalent polybox  codes $V,W\subset S^d$ it is possible to pass from $V$ to $W$ by gluing and cutting?}

\medskip
In the case of  cube tiling codes $V\subset S^d$ (that is $|V|=2^d$), where $S=\{a,a',b,b'\}$ the above problem was posed and resolved for $d\leq 4$ by Dutour, Itho and Poyarkov in \cite{DIP}. They computed that for every cube tiling codes $V,W\subset S^d$, $d\leq 4$, one can pass from $V$ to $W$ by gluing and cutting (these authors called such operation {\it flipping}). They asked also in \cite{DIP} whether the same is possible for $d=5,6,7$ and $S=\{a,a',b,b'\}$. In \cite{Ost}  Mathew, {\"O}sterg{\aa}rd and Popa gave an affirmative answer for $d=5$ (they called gluing and cutting {\it shifting}). In \cite{Kii} we showed that the above results are true for any $S$. On the other hand, for all $d\geq 8$ there are cube tiling codes $V,W\subset S^d$ such that $V$ cannot be obtained from $W$ by gluing and cutting (\cite{LS1,M}).

There are others tilings which can be modify by a sequence of local transformations. The most known are 2-dimensional domino and lozenge tilings (\cite{AS,R,W}). Local transformations of  domino tilings in dimensions $d\geq 3$ are also examined (\cite{KS,MS}).

\begin{uw}
{\rm We emphasize that passing from a twin pair $v,w$ with $v_i=w_i'\neq *$ to a twin pair $q,t$ with $q_i=t_i'\neq *$ and $v_j=q_j$ for all $j\neq i$, we may create the gluing $u$, because sometimes application of the gluing operation only (without cutting) may lead to interesting observations on the structure of polybox codes, as it was done in \cite{Kii}. However, in the presented paper we shall not transform proper polybox codes into improper (that is codes containing words $u\in (*S)^d$ with $u_i=*$ for some $i\in [d]$ ).}
\end{uw}

\begin{nap}{\rm 
Let
$$
V=\{aa, aa',a'b, a'b'\}\;\; {\rm and}\;\; W=\{cc, c'c, bc', b'c'\}.
$$
A realization of the cube tiling code $V$ is pictured in Figure 1 (on the left), and on the right a realization of $W$ is presented. Between them, we have realizations of three cube tiling codes: $V^1=V\setminus \{aa, aa'\} \cup \{ac, ac'\}$, $V^2=V^1\setminus \{a'b, a'b'\} \cup \{a'c, a'c'\}$ and $V^3=V^2\setminus \{ac', a'c'\} \cup \{bc', b'c'\}$. To pass from $V$ to $W$ by gluing and cutting we made the following sequence of local transformations:
$$
\{aa, aa'\}\rightarrow \{ac, ac'\},\;\; \{a'b, a'b'\}\rightarrow \{a'c, a'c'\},\;\; \{ac', a'c'\}\rightarrow \{bc', b'c'\}, \;\; \{ac, a'c\}\rightarrow \{cc, c'c\}.
$$ 
}
\end{nap}

\vspace{-0mm}
{\center
\includegraphics[width=12cm]{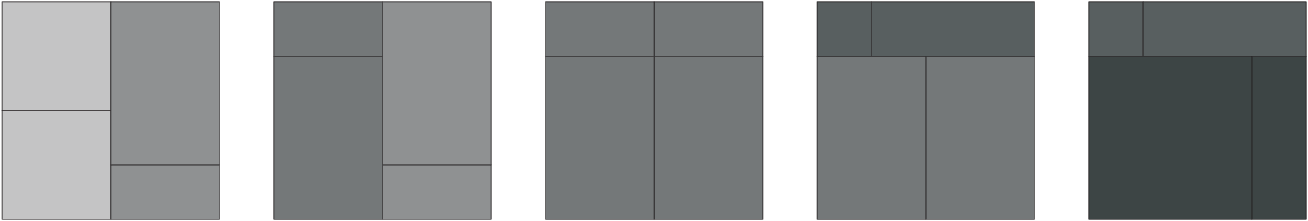}\\
}

\medskip
\noindent{\footnotesize Figure 1: Realizations of the codes (from the left to the right) $V,V^1,V^2,V^3$ and $W$.   In the realizations of the codes we took $X=[0,1]^2$, $f_1(a)=f_2(a)=[0,\frac{1}{2}), f_1(b)=f_2(b)=[0,\frac{1}{4})$ and $f_1(c)=f_2(c)=[0,\frac{3}{4})$.
}

\medskip
Clearly, if $W$ is obtained from $V$ by gluing and cutting, then both codes have to contain twin pairs. Thus, if $V$ or $W$ does not contain a twin pair, then none of them can be obtained from the other by gluing and cutting.  In Theorem \ref{p12} we give a complete description of twin pair free equivalent codes in dimensions up to six having at most 16 words which was obtained in \cite[Theorem 2.7]{Kap}, \cite[Theorem 31]{KisL} and \cite[Theorem 7.1]{Kis}. To do this, we need the following notations: If $A=\{i_1<...<i_k\}\subset [d]$ and $v\in S^d$, then $v_A=v_{i_1}\ldots v_{i_k}\in S^k$ and $V_A=\{v_A\colon v\in V\}$. To simplify the notation let $v_{i^c}=v_{\{i\}^c}$, that is, the word $v_{i^c}\in S^{d-1}$ arises from $v$ by skipping the letter $v_i$ in $v$. Moreover, $V_{i^c}=\{v_{i^c}:v\in V\}$.

\begin{tw}
\label{p12}
For $d\leq 6$ there is one, up to isomorphism, pair of twin pair free disjoint and equivalent polybox codes $V,W\subset S^d$ with $|V|\leq 16$. The codes $V,W$ have the form  

$$
V_A=\{aa'bb', abb'a, ab'b'b', a'ab'b', a'a'ab', a'bb'b, babb', bbbb, bb'a'b, b'aba', b'a'bb, b'b'b'b\}
$$
and
$$
W_A=\{a'a'a'b, a'baa', baa'a, aa'a'a, a'aa'a', abba', bbaa, ab'a'a', b'ab'a, b'a'aa, b'b'aa', bb'ab'\},
$$
where $A=\{1,2,3,4\}\subset [d]$ and $V_{A^c}=W_{A^c}=\{r_{A^c}\}$, $A^c=[d]\setminus A$, where $r\in S^d$ is any fixed word.
\hfill{$\square$}
\end{tw}

We shall call the pair $V,W$ described in the above theorem a {\it special pair}.

Let $V,W\subset S^d$ be polybox codes, and let $V\sim W$ if and only if $V$ and $W$ are equivalent. Clearly, $\sim$ is an equivalence relation. On the equivalence class $[V]_{\sim}=\{W\subset S^d\colon W\sim V\}$ we define a new relation $\approx$: $P\approx Q$ if and only if $P$ can be obtained from $Q$ by gluing and cutting. The relation $\approx$ is also an equivalence relation. If $P\approx Q$, then we say that $P$ and $Q$ are {\it strongly equivalent}. Let $\ka T(S^d)$ be the family of all cube tiling codes $V\subset S^d$. Note that, since every two cube tiling codes in $\ka T(S^d)$ are equivalent, we may consider the quotient set $\ka T(S^d)/_{\approx}$.
Thus, the question whether every pair of cube tiling codes are strongly equivalent is a question about the cardinality of the set $\ka T(S^d)/_{\approx}$. For example, we know that $|\ka T(S^d)/_{\approx}|\geq 2$ for every $d\geq 8$: If $M\subset S^8$ is Mackey's counterexample to Keller's conjecture given in \cite{M}, that is $M$ is a cube tiling code without twin pairs, then $M\approx V$ if and only if $V=M$. 

Alternatively, we may ask on a connectivity of the following graph: Let $G$ be a graph with $\ka T(S^d)$ as the set of vertices. Two vertices $V,W\in \ka T(S^d)$  are joined if $V$ and $W$ differs by one operation in the glue and cut procedure, that is there are two twin pairs $v,w\in V$ and $p,q\in W$ such that gluing of $v,w$ is the same as gluing of $p,q$ and $V\setminus\{v,w\}=W\setminus\{p,q\}$. Thus, if $V\approx W$ for every $V,W\in \ka T(S^d)$, then the graph $G$ is connected (compare \cite{DIP,Ost}). Of course, the connectedness of the graph $G$ means that starting from any $V\in \ka T(S^d)$ we can generate, by gluing and cutting, the entire  set $\ka T(S^d)$. (Let us mention that for $|S|=64$ we have $|\ka T(S^6)|>10^{84}$ (\cite[Subsection 3.4]{Kii}).)


\smallskip
In the presented paper we prove the following theorem:

\begin{tw}
\label{gc}
Every cube tiling codes $V,W\subset S^6$ are strongly equivalent.
\end{tw}


\medskip
\medskip
Let $V,W\subset S^d$ be polybox codes. We shall write $V\dot{\sqsubseteq} W$ if $V\sqsubseteq W$ and for every $v\in V$ it is not possible to pass from $W$ to some (equivalent) polybox code $\bar{W}$ by gluing and cutting, where $v\in \bar{W}$. Clearly, if $V\dot{\sqsubseteq} W$ and $W\dot{\sqsubseteq} V$, then the codes $V,W$ are not strongly equivalent.


Our proof of Theorem \ref{gc} is based on a partial characterization of strongly equivalent polybox codes $V,W\subset S^5$.  More precisely, we shall prove that:

\begin{tw}
\label{al}
Let $V,W\subset S^5$ be equivalent polybox codes such that $|V|\leq 15$, $V\dot{\sqsubseteq} W$ and $W\dot{\sqsubseteq} V$. 
\begin{enumerate}
\item If $|S|\geq 8$ and $|V|\leq 8$ or $|S|=6$ and $|V|\leq 10$, then there are no such codes $V,W$. In particular, equivalent polybox codes $V,W\subset S^5$, $|S|\leq 6$, with $|V|\leq 10$ are strongly equivalent. 
\item If $S=\{a,a',b,b'\}$, then there  is precisely one, up to isomorphism, such a pair of codes $V,W$ and they form the special pair. In particular, every equivalent codes $V,W\subset S^5$, $S=\{a,a',b,b'\}$, with $|V|\leq 11$ are strongly equivalent.
\end{enumerate}
\end{tw}

\begin{uw}{\rm
We believe that Theorem \ref {al} is true without any restriction on the cardinality of $S$, that is: {\it If $V,W\subset S^5$ are equivalent polybox codes such that $|V|\leq 11$, then $V$ and $W$ are strongly equivalent.} However, computations in this general case are longer, and to prove Theorem \ref{gc} we need Theorem \ref{al} in the form given above. }
\end{uw}

In $1930$, Keller conjectured that in every cube tiling of $\er^d$ there is a twin pair (\cite{Ke1}). It was known that Keller's conjecture is true for dimensions $d\leq 6$ (\cite{P}) and false for all dimensions $d\geq 8$ (\cite{LS1,M}). Recently, Brakensiek, Heule, Mackey and Narv{\'a}ez made in outstanding way the final step in proving Keller's conjecture in dimension seven (\cite{BHM}).
In the three papers on Keller's conjecture in dimension seven \cite{Kap,Kis,KisL} and also in \cite{KP}, we developed a method of analysing the structure of polybox codes via examination of covers of a polybox code by another polybox code. Some aspects of this `covering method' will be used in the paper. Let us recall that the method applied in \cite{Ost} for cube tiling codes $V\subset S^5$, where $S=\{a,a',b,b'\}$ was based on enumeration of all non-isomorphic such codes (there are $899,710,227$ of them (\cite{Ost})). It seems that an attempt to prove Theorem \ref{gc} via enumeration of all non-isomorphic cube tiling codes $V\subset S^6$, $|S|\leq 64$, is doomed to failure: The number of all cube tiling codes $V\subset S^5$, where $|S|=4$ is of the order $6\cdot 10^{14}$ (\cite{Ost}), while the number of all cube tiling codes $V\subset S^6$, where $|S|=64$, as we mentioned above, far exceeds the number $10^{84}$.

At the end of this section, let us note that the problem of passing from a cube tiling code $V$ to a given cube tiling code $W$ by gluing and cutting is very similar (taking into account a manner of doing a single operation) to the popular 15-puzzle game (\cite{Wo}). (Recently, a three dimensional version of this game, called Varikon cube was examined by d'Eon and Nehaniv in \cite{EN}.) In our {\it glue and cut game} the object is to reach a given cube tiling code $V\subset S^d$ by gluing and cutting starting from some other cube tiling code $W\subset S^d$. (Clearly, we may consider realizations of cube tiling codes (minimal partitions, in particular cube tilings of $\te^d$), such as in Figure 1 and 2, which allow us to make transformations of one twin pair into the another in a continuous fashion.)  
Similarly like in 15-puzzle, there are unsolved configurations in the glue  and cut game. As we mentioned, if $M\subset S^8$ is Mackey's counterexample to Keller's conjecture, then for every $W\subset S^8$, $W\neq M$, $W$ and $M$ are not strongly equivalent. On the other hand the results in \cite{DIP,Ost} and Theorem \ref{gc} show that, every two cube tiling codes $V,W\subset S^d$, $d\leq 6$, form a solved configuration in the glue and cut game.

\section{Basic notions on polybox codes}

As it was mentioned in the previous section, we can interpret a polybox code as a system of boxes. There are many such interpretations (realizations), but we shall use the following one which has particular nice properties (\cite[Section 10]{KP}). Let
$$
 ES=\{B\subset S\colon |\{s,s'\}\cap B|=1, \text{whenever $s\in S$}\}\;\;\;{\rm and}\;\;\; E s=\{B\in ES\colon s\in B\}.
$$
Let $V\subseteq S^{d}$ be a polybox code, and let $v\in V$. The {\it equicomplementary realization} of the word $v$ is the box   
$$
\breve{v}=Ev_1\times \cdots \times Ev_d
$$
in the $d$-box $(ES)^d =ES\times \cdots \times ES.$ The equicomplementary realization of the code $V$ is the family
$$
E(V)=\{\breve{v}:v\in V\}.
$$ 
If $s_1,\ldots , s_n\in S$ and $s_i\not\in\{s_j, s'_j\}$ for every $i\neq j$, then
\begin{equation}
\label{dkostki}
|E s_1 \cap\dots\cap E s_n|=(1/2^{n})|ES|.
\end{equation}
The value of the realization $E(V)$, where $V\subseteq S^d$,  lies in the equality (\ref{dkostki}). 
In particular, boxes in $E(V)$ are of the same size: $|E v_i|=(1/2)|ES|$ for every $i\in [d]$ and consequently  $|\breve{v}|=(1/2^d)|ES|^d$ for $v\in E(V)$. Thus, two boxes $\breve v, \breve w\subset (ES)^d$ are dichotomous, and also the words $v,w$ are dichotomous, if and only if $\breve v\cap \breve w=\emptyset.$ The same is true for translates of the unit cube in the flat torus $\te^d$ and therefore working with the boxes $\breve{v}, v\in V$,  we can think of them as translates of the unit cube in $\te^d$.

If $w\in S^d$ and  $V\subset S^d$ is a polybox code, then it can be shown (\cite[Theorem 10.4]{KP}) that $w\sqsubseteq V$ if and only if $\breve{w}\subseteq \bigcup E(V)$. A cover $V$ of $w$ is {\it minimal} if $\breve{w}\cap \breve{v}\neq\emptyset$ for every $v\in V$. Similarly, a polybox code $C_P\subset S^d$  is a \textit{minimal} cover of a code $P\subset S^d$ if $P\sqsubseteq C_P$ and for every $c\in C_P$ there is $p\in P$ such that $\breve{c}\cap \breve{p}\neq\emptyset$.

The following characterization of covers of words will be very useful in our computations (Subsection 3.1)

Let  
$g\colon S^d\times S^d\to \zet$ be defined by the formula $g(v,w)=\prod^d_{i=1}(2[v_i=w_i]+[w_i\not\in\{v_i, v'_i\}])$,
where $[\cdot]$ denotes the {\it Iverson bracket}, that is $[p]=1$ if the sentence $p$ is true and $[p]=0$ if it is false. Let 
\begin{equation}
\label{pv}
|w|_V=\sum_{v\in V} g(v,w).
\end{equation}
In \cite[Theorem 10.4]{KP} it was showed that  
\begin{equation}
\label{2d}
w \sqsubseteq V \Leftrightarrow |w|_V= 2^d.
\end{equation}


\subsection{Isomorphisms}

If $v\in (*S)^d$, and $\sigma$ is a permutation of the set $[d]$, then $\bar{\sigma}(v)=v_{\sigma(1)}\ldots v_{\sigma(d)}$. For every $i\in [d]$ let $h_i:*S\rightarrow *S$ be a bijection such that $h_i(l')=(h_i(l))'$ for every $l\in *S$ and $h_i(*)=*$. 
Let $h:(*S)^d\rightarrow (*S)^d$ be defined by the formula $h(v)=h_1(v_1)\ldots h_d(v_d)$. 
The group of all possible mappings $h\circ \bar{\sigma}$  will be denoted by $G((*S)^d)$ or $G(S^d)$ depending on whether we consider words written down in the alphabet $*S$ or in $S$. Let $S$ and $T$ be two alphabets, and let $|S|\leq |T|$. Two polybox codes $V\subset (*S)^d$ and $U\subset (*T)^d$ are {\it isomorphic} if there is  $h\circ \bar{\sigma}\in G((*T)^d)$ such that $U=h\circ \bar{\sigma}(\tau^d(V))$, where $\tau\colon *S\rightarrow *T$ is a fixed injection such that $\tau(*)=*$, $\tau(s')=\tau(s)'$ for $s \in S$ and $\tau^d(v)=\tau(v_1)\ldots \tau(v_d)$ for $v\in (*S)^d$. The composition $h\circ \bar{\sigma}$ is an {\it isomorphism} between $V$ and $U$. 

It is easy to check that for every polybox codes $V,W\subset (*S)^d$ and every $g\in G((*S)^d)$ if $V\sqsubseteq W$, then $g(V)\sqsubseteq g(W)$ (compare \cite[Section 5]{KisL}). Thus, if codes $V,W$ are equivalent and $g\in G((*S)^d)$, then the codes $g(V)$ and $g(W)$ are equivalent too. Similarly, if $V$ and $W$ are strongly equivalent, then $g(V)$ and $g(U)$ are strongly equivalent. 

\subsection{Distribution of words and passing from a cube tiling code to a cube tiling code}

If $V\subseteq S^d$, $l\in S$ and $i\in [d]$, then $V^{i,l}=\{v\in V\colon v_i=l\}$. If $S=\{a_1,a'_1,...,a_k,a_k'\}$, then the representation 
$$
V=V^{i,a_1}\cup V^{i,a_1'}\cup \cdots \cup V^{i,a_{k}}\cup V^{i,a_{k}'}
$$
will be called a {\it distribution of words in} $V$ (clearly, some sets $V^{i,a_j}$, $j\in [k]$, can be empty). If $R\subset S^d$ is a code such that for every $i\in [d]$ there is $j\in [k]$ such that $R=R^{i,a_j}\cup R^{i,a'_j}$, then $R$ is called a {\it simple code}. If $R$ is a simple cube tiling code and $V\cap R\neq\emptyset$, then the code $V\cap R$ is called a {\it simple component} of $V$ (compare Figure 2).

It is rather obvious that if $V$ and $W$ are cube tiling codes which are simple, then one can pass from $V$ to $W$ by gluing and cutting. Thus, to show that cube tiling code $V$ can be  transformed to a cube tiling code $W$ by gluing and cutting it is enough to show that both codes can be converted by gluing and cutting into simple codes.   Now we present the structure of cube tiling codes which explains our interest in Question 1 in the light of gluing and cutting of cube tilings codes.

\smallskip
Let us recall that $v_{i^c}=v_1\ldots v_{i-1}v_{i+1}\ldots v_d$ for $v\in S^d$ and $V_{i^c}=\{v_{i^c}:v\in V\}\subset S^{d-1}$ for $V\subset S^d$.

\smallskip
Let $U=U^{i,a_1}\cup U^{i,a_1'}\cup \cdots \cup U^{i,a_{k}}\cup U^{i,a_{k}'}$ be a cube tiling code.
It is known (\cite[Section 2]{Kii}) that for every $i\in [d]$ and every $j\in [k]$ the polybox codes $U^{i,a_j}_{i^c},U^{i,a_j'}_{i^c}\subset S^{d-1}$ are equivalent. Let $V=U^{i,a_k}_{i^c}$ and $W=U^{i,a_k'}_{i^c}$ and assume that $V\neq\emptyset$. (Later in the text, if we do some operation on such defined sets $V$ and $W$, then we shall  assume tacitly that they are non-empty sets.)
Since $V,W\subset S^{d-1}$ are equivalent, we may try to pass from $V$ to $W$ by gluing and cutting. Suppose that such passing is possible. 
Clearly, this means that we may pass from $U$ to a cube tiling code $P$ by gluing and cutting, where $P$ is such that 
$$
P=U^{i,a_1}\cup U^{i,a_1'}\cup \cdots \cup U^{i,a_{k-1}}\cup U^{i,a_{k-1}'}\cup P^{i,a_{k}}\cup U^{i,a_{k}'}\;\;\; {\rm and}\;\;\; P^{i,a_{k}}_{i^c}=U^{i,a'_{k}}_{i^c}.
$$
The last equality means that the code $P^{i,a_{k}}\cup U^{i,a_{k}'}$ consists of twin pairs $p,u$, $p\in P^{i,a_{k}}, u\in U^{i,a_{k}'}$ in the $i$th direction (Figure 2, the third picture from the left in the lower row). Now each pair $p,u$ can be glued and cut into the twin pair $\bar{p},\bar{u}$ with $\bar{p}_i=a_{k-1}$ and $\bar{u}_i=a'_{k-1}$ (Figure 2, the last three pictures in the lower row). This means that we may pass, by gluing and cutting,  from $P$ to a cube tiling code $\bar{U}$, where $\bar{U}$ has the following distribution:
$$
\bar{U}=U^{i,a_1}\cup U^{i,a_1'}\cup \cdots \cup U^{i,a_{k-2}}\cup U^{i,a_{k-2}'}\cup \bar{U}^{i,a_{k-1}}\cup \bar{U}^{i,a_{k-1}'}.
$$

Since the operation of gluing and cutting is transitive, we may pass from $U$ to $\bar{U}$ by gluing and cutting. Thus, we passed from $U$ to a simpler cube tiling code $\bar{U}$ in which the number of components in the distribution of words in $\bar{U}$ in the $i$th position is lower than that in $U$ (Figure 2). 

\vspace{-0mm}
{\center
\includegraphics[width=12cm]{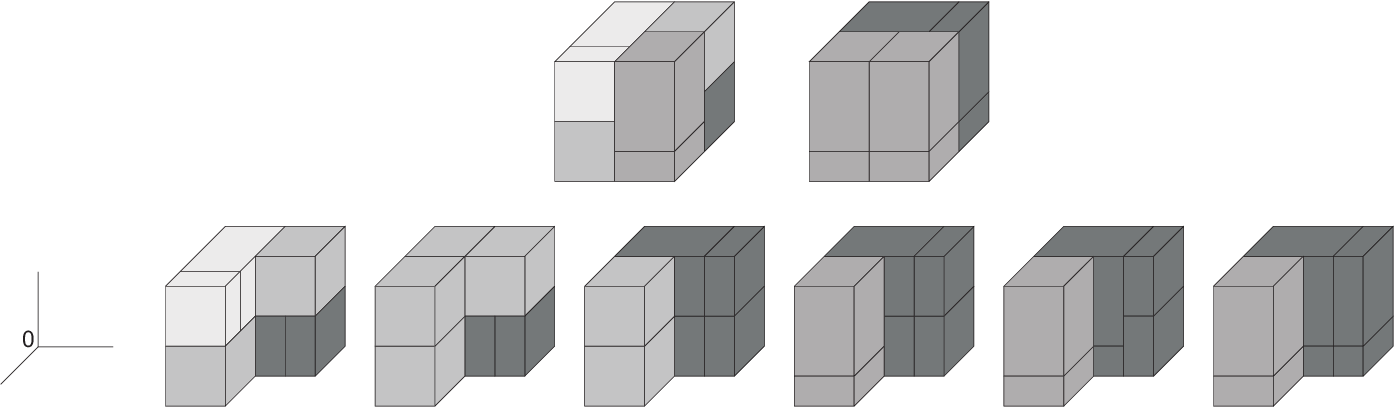}\\
}

\medskip
\noindent{\footnotesize Figure 2:   A passing from the cube tiling code $U=U^{3,a}\cup U^{3,a'}\cup U^{3,b}\cup U^{3,b'}$ (the picture on the left in the upper row)  to $\bar{U}=\bar{U}^{3,a}\cup \bar{U}^{3,a'}$ (the picture on the right in the upper row), where $U=\{aab, a'a'b', ab'b', abb', ba'b, b'a'b, a'aa, a'aa'\}$ and $\bar{U}=\{aaa, aaa',a'aa, a'aa', ba'a,\\ ba'a', b'a'a, b'a'a'\}$. Simple components of realizations of $U$ and $\bar{U}$ are one-colored. In the first three pictures in the lower row we present a passing by gluing and cutting from $ U^{3,b}\cup U^{3,b'}$ to $P^{3,b}\cup U^{3,b'}$, and the last three pictures in the lower row is a passing from mentioned above twin pairs $p,u$ to $\bar{p},\bar{u}$. In the realizations we took $X=[0,1]^3, f_1(a)=f_2(a)=[0,\frac{1}{2}),f_3(a)=[0,\frac{1}{4})$ and $f_1(b)=f_2(b)=[0,\frac{3}{4}), f_3(b)=[0,\frac{1}{2})$.
}

\medskip

Along this lines we could try to reduce by gluing and cutting the code $U$ to some simple cube tiling code. But this depends on our knowledge about strongly equivalent polybox codes $V,W\subset S^{d-1}$.
Theorem \ref{al} provides us with such knowledge.

In the next section we list results which are needed in our proof of Theorem \ref{al}.


\section{Preliminary results}

Our goal is to characterize the family of strongly equivalent codes $V,W\subset S^5$ with $|V|\leq 15$, but we start with a brief explanation of how we can construct a pair of equivalent codes $V,W$ having some property, say a property $\bf{P}$. The first step is to establish initial codes $V_0\subset V$ and $W_0\subset W$, as large as possible, which have to appear in $V$ and $W$, and next based on the codes $V_0$ and $W_0$ try to construct $V$ and $W$. The simplest assumption is $V_0=\{v\}$, where $v$ is a word. Now we can say quite a lot about $W_0$: We can assume that $W_0$ is a minimal cover of $v$. If for example the property $\bf{P}$ means that $V$ and $W$ are twin pair free, then the cover $W_0$ must be twin pair free. It is known that in such a case $|W_0|\geq 5$, and the family of all non-isomorphic twin pair free minimal covers of $v$ for $d\leq 5$ and $|W_0|\leq 15$ can be easily computed (Statement \ref{sl11}). 

Of course, strongly equivalent codes $V,W$ contain twin pairs but as we shall show below we  may assume that some words in $V$ are covered by codes without twin pairs. This reduces a number of initial configurations which have to be considered and makes computations possible in a reasonable time.



\begin{lemat}
\label{bez}
Let $P\subset S^d$ be a simple polybox code, and let $C_P\subset S^d$ be a cover of $P$. Then one can pass by gluing and cutting from $C_P$ to a cover $\bar{C}_P$ of $P$  such that for every $p\in P$ the minimal cover $\bar{C}_p=\{w\in \bar{C}_P\colon \breve{w}\cap \breve{p}\neq\emptyset\}$ of $p$ does not contain a twin pair.
\end{lemat}
\proof For simplicity of the notation we may assume that $P\subset \{a,a'\}^d$. Let $p\in P$, and let $v,w\in C_P$ be a twin pair with $v_i'=w_i$ such that $\breve{p}\cap \breve{v}\neq\emptyset$ and $\breve{p}\cap \breve{w}\neq\emptyset$. Clearly, $p_i\not\in \{v_i,w_i\}$. Let us pass from $\{v,w\}$ to the twin pair $\{r,q\}$ by gluing and cutting such that $r_i=a,q_i=a'$ (and of course $r_j=v_j$ for every $j\neq i$), and let $C^1_P=C_P\setminus \{v,w\}\cup \{r,q\}$. For every $t\in P$ the twin pair $r,q$ is not contained in the minimal cover $C_t\subset  C^1_P$ of $t$ (because $r_i=a,q_i=a'$). 
If for some $t\in P$ the minimal cover $C_t\subset C^1_P$ of $t$ contains a twin pair we pass in the similar way from $C^1_P$ to a code $C^2_P$. Since during this process we always change a pair of letters $s,s'\in S$ into the pair $a,a'$ (in transformed twin pair), after $n$ steps we have to obtain a cover $C^n_P$ of $P$ such that for every $p\in P$ the minimal cover $C_p\subset C^n_P$ of $p$ does not contain a twin pair. Then $\bar{C}_P=C^n_P$.
\hfill{$\square$} 

\medskip
In \cite[Theorem 10.6]{KP} we proved the following structural result which is a base in the construction of our initial configurations $V_0,W_0$.

\begin{lemat}
\label{ft}
Let $C_p\subset S^d$ be a twin pair free cover of a word $p\in S^d$. There are $v,w\in C_p$ such that $v_i=w_i$ or $v_i=w'_i$ for every $i\in [d]$ and  $|\{i\in [d]\colon v_i=w_i'\}|=2k+1$ for some positive integer $k\geq 1$.
\hfill{$\square$}
\end{lemat}

\medskip
Let $V,W\subset S^d$ be polybox codes. Recall that if $V\dot{\sqsubseteq} W$ and $W\dot{\sqsubseteq} V$, then the codes $V,W$ are not strongly equivalent. Note that, for equivalent codes $V,W$ it can happen that $V\dot{\sqsubseteq} W$ but not conversely. For example, if $M,U\subset S^8$ are cube tiling codes, where $M$ is Mackey's counterexample (\cite{M}), and $U$ is any simple cube tiling code, then we have only $U\dot{\sqsubseteq} M$. However, if equivalent codes $V,W\subset S^d$ are not strongly equivalent, then they have to be related to some codes $V^1,W^1$ with $V^1\dot{\sqsubseteq} W^1$ and $W^1\dot{\sqsubseteq} V^1$. More precisely, we have the following lemma:

\begin{lemat}
\label{sind}
If equivalent polybox codes $V,W\subset S^d$ are not strongly equivalent, then it is possible to pass by gluing and cutting from $V$ to $\bar{V}$ and similarly from $W$ to $\bar{W}$ such that there are decompositions $\bar{V}=V^1\cup P$ and $\bar{W}=W^1\cup P$, where $V^1\neq\emptyset$, $V^1\dot{\sqsubseteq} W^1$ and $W^1\dot{\sqsubseteq} V^1$.
\end{lemat}
\proof
If $V\dot{\sqsubseteq} W$ and  $W\dot{\sqsubseteq} V$, then $V^1=V,W^1=W$ and $P=\emptyset$. If it is not true that $V\dot{\sqsubseteq} W$, then there is  $v\in V$ such that $v\in R$, where $R$ is obtained from $W$ by gluing and cutting and the equivalent codes $V_1=V\setminus \{v\}$, $W_1=R\setminus \{v\}$ are not strongly equivalent. By induction on the number of words in $V$, it is possible to pass by gluing and cutting from $V_1$ and $W_1$ to $\bar{V}_1$ and $\bar{W}_1$, respectively and $\bar{V}_1=V^1\cup Q$ and $\bar{W}_1=W^1\cup Q$, where $V^1\neq\emptyset$, $V^1\dot{\sqsubseteq} W^1$, $W^1\dot{\sqsubseteq} V^1$. Thus, $\bar{V}=V^1\cup P$ and $\bar{W}=W^1\cup P$, where $P=Q\cup \{v\}$.
\hfill{$\square$}

\medskip
In the next lemma we establish our initial configurations.

\begin{lemat}
\label{pbez}
Let $U,R\subset S^5$ be two equivalent polybox codes such that $U\dot{\sqsubseteq} R$ and $R\dot{\sqsubseteq} U$, and let $P=\{bbbbb, \\b'b'b'bb\}$ or $P=\{bbbbb,b'b'b'b'b'\}$. There is a mapping  $g\in G(S^d)$ such that the equivalent codes $V=g(R)$ and $W=\overline{g(U)}$, where $\overline{g(U)}$ is obtained from $g(U)$ by gluing and cutting, have the following properties: $P\subset V$, the minimal cover $C_P\subset W$ of $P$ is such that the minimal covers $C_s\subset C_P$ of $s\in \{bbbbb,b'b'b'bb, b'b'b'b'b'\}$ do not contain a twin pair. Moreover, $V\dot{\sqsubseteq} W$, $W\dot{\sqsubseteq} V$.

\end{lemat}
\proof
Fix any $u\in U$ and let $C_u\subset R$ be the minimal cover of $u$. There are two possibilities: First, $C_u$ does not contain a twin pair. Then, by Lemma \ref{ft}, there is a polybox code $Q\subset C_u$ and a mapping $g\in G(S^5)$ such that $g(Q)=P$. Second, $C_u$ contains a twin pair. Then, by Lemma \ref{bez}, we pass by gluing and cutting from $C_u$ to a cover (not minimal) $\bar{C}_u$ of $u$ such that the minimal cover $M_u\subset \bar{C}_u$ of $u$ has more than one word (we can do this as $U\dot{\sqsubseteq} R$ and $R\dot{\sqsubseteq} U$) and does not contain a twin pair. Consequently, we may assume that $Q\subset M_u$. (Clearly, passing from $C_u$ to $\bar{C}_u$ means that we pass from $R$ to $\bar{R}$ by gluing and cutting and $\bar{C}_u\subset \bar{R}$.) Let $U_1=g(U)$ and $R_1=g(R)$ (in the first case)  or $U_1=g(U)$ and $R_1=g(\bar{R})$ (in the second case). We have $P\subset R_1$ and $U_1,R_1$ are equivalent (Subsection 2.1). Now we take the minimal  cover $C_P\subset U_1$ of $P$. If covers $C_s\subset C_P$ of words $s\in P$ do not contain a twin pair, then $V=R_1$ and $W=U_1$. If it is not so, by Lemma \ref{bez}, we pass by gluing and cutting from $C_P$ to $\bar{C_P}$ such that covers $\bar{C}_s\subset \bar{C}_P$ of words $s\in P$ do not contain a twin pair. This passing means that we pass from $U_1$ to a code $\bar{U_1}$ by gluing and cutting. Then we take $V=R_1$ and $W=\bar{U_1}$, and consequently $P\subset V$ and $\bar{C}_P\subset W$.

Obviously, if $\bar{R}$ is obtained from $R$ by gluing and cutting, then, since $U\dot{\sqsubseteq} R$ and $R\dot{\sqsubseteq} U$, we have $U\dot{\sqsubseteq} \bar{R}$ and $\bar{R}\dot{\sqsubseteq} U$. Moreover, any $g\in G(S^d)$ preserves equivalency of polybox codes (compare Subsection 2.1). Thus, $V\dot{\sqsubseteq} W$ and $W\dot{\sqsubseteq} V$.
\hfill{$\square$}

\medskip
Let $v=bbbbb$, $w=b'b'b'bb$, $u=b'b'b'b'b'$, and let $P=\{v,w\}$, $T=\{v,u\}$. The above results show that looking for strongly equivalent codes $V,W\subset S^5$ we may start with initial configurations $V_0,W_0$ of the form: $V_0=P$ or $V_0=T$, and $W_0$ is a cover of $V_0$ such that minimal covers $C_v,C_w$, $C_u$ do not contain twin pairs, where $C_v,C_w\subset W_0$ if $V_0=P$ and $C_v,C_u\subset W_0$ if $V_0=T$. The following two observations allow us to make one more useful assumption on $V,W$:

\begin{lemat}
\label{four}
Let $v\in S^d$, and let $C_v\subset S^d$ be a minimal cover of $v$ having less than five words. Then we can pass from $C_v$ to a polybox code $\bar{C_v}$ by gluing and cutting, where  $\bar{C_v}$ contains $v$. 
\end{lemat}
\proof 
It is an immediate consequence of the fact that  every cover of $v$ having less than five words contains a twin pair (\cite[Corollary 3.3]{Kap}). Thus, if $w,u\in C_v$ is a twin pair with  $w_i=u_i'$, $w_i,u_i\neq v_i$, then passing from  $\{w,u\}$ to $\{\bar{w},\bar{u}\}$, where $\bar{w}_i=v_i$ and $\bar{u}_i=v_i'$, by gluing and cutting we get the cover $C^1_v=C_v\setminus \{w,u\} \cup \{\bar{w},\bar{u}\}$ of $v$ in which $|\{q\in C^1_v\colon \breve{v}\cap \breve{q}\neq\emptyset\}|=|C_v|-1$ (as $\breve{v}\cap \breve{\bar{u}}=\emptyset$). In this manner we can obtain $\bar{C}_v$.
\hfill{$\square$}

\medskip
Let $V,W\subset S^d$ be polybox codes.  The number
$$
\delta(V,W)=\min_{w\in W}|\{v\in V\colon v\; {\rm and}\; w\; {\rm are \; not\; dichotomous}\}| 
$$
will be called a {\it density of $V$ with respect to $W$}. Note that the density has a clear geometrical meaning: $\delta(V,W)=\min_{w\in W}|\{v\in V\colon \breve{w}\cap \breve{v}\neq\emptyset\}|$ (compare the beginning of Section 2). The number  $\Delta(V,W)=\min\{\delta(V,W),\delta(W,V)\}$ will be called a {\it common density} of the codes $V$ and $W$. 

Thus, from Lemma \ref{four} we get

\begin{wn} 
\label{del4}
If $V,W\subset S^d$ are polybox codes, $V$ is a cover of $W$ and $\delta(V,W)\leq 4$, then there is $w\in W$ such that one can pass by gluing and cutting from $V$ to $\bar{V}$, where $\bar{V}$ contains $w$.\hfill{$\square$}
\end{wn}

In the next result we give the structure of polybox codes which allows us to simplify some computations. If $V\subset S^d$ and $V=V^{i,s}$ for some $i\in [d]$ and $s\in S$, then $V$ is called {\it flat}. From the definition of equivalency of codes and (\ref{dkostki}), it follows that if $V$ is a flat code, and $W$ is an equivalent code  to $V$, then $W$ is also flat (that is $W=W^{i,s}$).

\begin{lemat}
\label{fla}
Let $S=\{a,a',b,b'\}$, and let $V\subset S^d$ be a polybox code. If $V=V^{i,a}\cup V^{i,a'}\cup V^{i,b}\cup V^{i,b'}$ for some $i\in [d]$,  $1\leq |V^{i,a'}|+|V^{i,b}|+|V^{i,b'}|\leq 4$  and  $W$ is an equivalent code to $V$,  then $\Delta(V,\bar{W})\leq 4$, where $\bar{W}$ is a code obtained from $W$ by gluing and cutting.
\end{lemat}
\proof
Assume first that $W^{i,a'}=\emptyset$. Then $W^{i,b}\cup W^{i,b'}\neq\emptyset$, otherwise $W$ is flat, and $V$ is not, which is impossible.

If $V^{i,b}\cup V^{i,b'}=\emptyset$, then the codes $W^{i,b}_{i^c}$ and $W^{i,b'}_{i^c}$ are equivalent (\cite[Subsection 2.8, \bf{C}]{KisL}). Since $W^{i,a'}=\emptyset$, we have $|V^{i,a'}|+|V^{i,b}|=|W^{i,b}|$, and thus $|W^{i,b}|\leq 4$. By Lemma \ref{four}, we can pass from $W$ to a code $\bar{W}$ by gluing and cutting such that $\bar{W}=\bar{W}^{i,a}\cup \bar{W}^{i,a'}$ (compare Subsection 2.2). Thus, $V^{i,a'}$ and $\bar{W}^{i,a'}$ are equivalent, and since $1\leq |V^{i,a'}|\leq 4$, we have $\Delta(V,\bar{W})\leq 4$.

Therefore, we assume that $V^{i,b}\neq\emptyset$.  We have $v\sqsubseteq W^{i,b}$ for any $v\in V^{i,b}$, and then $\Delta(V,W)\leq 4$, as $|V^{i,a'}|+|V^{i,b}|\leq 4$ and $|V^{i,a'}|+|V^{i,b}|=|W^{i,b}|$.

Let now $W^{i,a'}\neq\emptyset$. For $w \in W^{i,a'}$, we have $w\sqsubseteq V^{i,a'}\cup V^{i,b}\cup V^{i,b'}$, and thus $\Delta(V,W)\leq 4$.
\hfill{$\square$}

\medskip
The next lemma is rather obvious.

\begin{lemat}
\label{ind}
Let $V,W\subset S^d$ be two equivalent polybox codes with $n$ words. Suppose that every two equivalent codes with $n-1$ words are strongly equivalent. If there is $w\in W$ such that $w\in \bar{V}$, where $\bar{V}$ is obtained from $V$ by gluing and cutting, then $V$ and  $W$ are strongly equivalent.\hfill{$\square$}
\end{lemat}

\smallskip
If we want to construct a polybox $V$ which is equivalent to a given code $W$, it is useful to apply the following observation expressed in Lemma \ref{index} (compare algorithm {\sc FindSecondCode} in the next subsection). 

Let $b\colon S\rightarrow \{0,1\}$ be such that $b(s)+b(s')=1$ for every $s\in S$. A {\it binary code} of a word $v\in S^d$ is the vector  $\beta(v)=(b(v_1),...,b(v_d))$. If $V\subset S^d$ is a polybox code, then $\beta(V)=\{\beta(v)\colon v\in V\}$. In \cite[Section 5]{KP} we proved

\begin{lemat}
\label{index}
If $V,W\subset S^d$ are equivalent polybox codes, then $\beta(V)=\beta(W)$.
\hfill{$\square$}
\end{lemat}

\medskip
Our proof of Theorem \ref{al} needs computer support. In the next subsection we list algorithms which shall be used.

\subsection{Algorithms}

Since we are interested in twin pair free covers of a word, in the following algorithm we show how to find such covers. 
The algorithm given below is closely related to the relation (\ref{2d}). (A deeper explanation can be found at the end of Subsection 2.8 in \cite{KisL}.)

\medskip
{\bf Algorithm} {\small {\sc CoverWord}}.

\medskip
Let $u\in S^5$, $k\geq 5$ be an integer, and let $\ka C^k$ be the family of all $k$-elements twin pair free minimal covers $C_u$ of $u$, that is  $\breve{u}\cap \breve{v}\neq\emptyset$ for every $v\in C_u$ and every $C_u\in \ka C^k$. We may assume that $u=bbbbb$. Moreover, by Lemma \ref{ft}, we assume that codes in $\ka C^k$ contain one of the codes: 
$$
V_{3,0}=\{aaaaa, a'a'a'aa\}, V_{3,1}=\{aaaab, a'a'a'ab\},
$$
$$
V_{3,2}=\{aaabb, a'a'a'bb\} \; {\rm or}\; V_{5,0}=\{aaaaaa, a'a'a'a'a'\}.
$$
Our goal is to find the family $\ka C^k$.

\smallskip
{\it Input}. The word $bbbbb\in S^5$ and the number $k$.

{\it Output}. The family $\ka C^k$.

\smallskip
1. Let $\ka S_k=\{(x_0,...,x_{4}) \in \en^5\colon  \sum_{i=0}^{4}x_i2^i=2^5\; {\rm and}\; \sum_{i=0}^{4}x_i=k\}$, where $\en=\{0,1,2,...\}$.

2. For $i \in \{0,...,4\}$ indicate the set $\ka A_i$ consisting of all words $v\in S^5$ such that $v$ contains precisely $i$ letters $b$ and $v$ does not contain the letter $b'$.

3. Fix $x\in \ka S_k$ and let $s(x)=\{i_1<\cdots <i_m\}$ consists of all $i_j\in \{0,...,4\}, j\in [m]$, for which $x_{i_j}>0$. Fix $V_{n,i_1}\subset \ka A_{i_1}$. 
For $i\in s(x)$ let $\ka B_i=\{v\in \ka A_i\colon V_{n,i_1}\cup \{v\}\;\; {\rm is\; a \; twin\;\; pair\;\; free\;\; code}\}$.

4. Let $I$ be the multiset containing $i_1$ with the multiplicity $x_{i_1}-2$ (recall that $x_{i_1}\geq 2$, as $x\in \ka S_k$) and  $i_j$ with the multiplicity $x_{i_j}$ for $j\in \{2,...,m\}$. By $I[j]$ we denote the $j$th element of $I$. 

5. Let $\ka D^2=\{V_{n,i_1}\}$.  

6. For $l\in \{2,...,k-1\}$ having computed $\ka D^l$ we compute the set $\ka D^{l+1}$: For $v\in \ka B_{I[l-1]}$ and for $U\in \ka D^l$ if $U\cup \{v\}$ is a twin pair free code, then we attach it to  $\ka D^{l+1}$. 

7. Clearly, $\ka D^k=\ka D^k(V_{n,i_1},x)$ so let $\ka C^k$ be the union of the sets $\ka D^k(V_{n,i_1},x)$ over $x\in \ka S^k$ and  $V_{n,i_1}\subset \ka A_{i_1}$ (recall that $V_{n,i_1}$ depends on $x$).


\bigskip
The next algorithm allows to check whether $\{p\}\dot{\sqsubseteq} V$. (Thus, it can be used to decide whether $V\dot{\sqsubseteq} W$.) Let us observe that the fact that a twin pair $\{r,q\}$ is obtained from a twin pair $\{v,w\}$ by gluing and cutting can be denoted by $\{v,w\} \sim \{r,q\}$ (that is, $\{v,w\}$ is equivalent to $\{r,q\}$). 

\medskip

{\bf Algorithm} {\small {\sc GlueAndCut}}.

\medskip
Let $\ka V$ be a family of polybox codes $V\subset S^d$, and let
$$
t(\ka V)=\bigcup_{V\in \ka V}\bigcup_{\substack{\{v,w\}\subset V \\ \{v,w\}\;{\rm is\;a \;twin\;pair}}}\{V\setminus\{v,w\}\cup \{r,q\}\colon \{r,q\}\subset S^d\;{\rm and}\;\{v,w\} \sim \{r,q\}\}.
$$
Moreover, let $t^0(\ka V)=t(\ka V)$ and $t^n(\ka V)=t(t^{n-1}(\ka V))$ for $n\geq 1$. Let $M$ be the smallest positive integer such that $t^n(\ka V)=t^{n+1}(\ka V)$ for $n\geq M$ (since $S$ is finite, such $M$ exists). Clearly, if $V\subset S^d$ is a polybox code and $\ka V=\{V\}$, then the set $t^M(\ka V)$ consists of all  polybox codes that can be obtained from $V$ by gluing and cutting (compare Section 1).

\medskip   
{\it Input}. A polybox code $V\subset S^d$ and a word $p\in S^d$ with $p\sqsubseteq V$.

{\it Output}. 1 if $\{p\}\dot{\sqsubseteq} V$, and 0 otherwise.

\medskip
1. For $n\in \{0,...,M\}$ and for $\bar{V}\in t^n(\{V\})$ if $\delta(\bar{V},\{p\})<5$, then return 0.

2. If $p\not\in \bar{V}$ for every $\bar{V}\in t^M(\{V\})$, then return 1.

\begin{uwi}
\label{u2}
{\rm Let us note that, by Lemma \ref{four}, it is better to check the condition $\delta(\bar{V},\{p\})<5$ rather than $p\in \bar{V}$. (Clearly, if $p\sqsubseteq \bar{V}$ and $\delta(\bar{V},\{p\})=1$, then $p\in \bar{V}$.)

The above algorithm is slow, but for small polybox codes $V$ it works quite well. However, when a polybox code $V$ contains a larger number of words its efficiency  goes down radically, as the number $m=|t^M(\ka V)|$ can be huge. In particular, for cube tiling codes and dimensions $d\geq 4$ it is rather worthless. For example, if $V\subset S^4$ is a cube tiling code and $|S|=16$, then  $m>10^{16}$ (\cite{Kii}), and if $V\subset S^5$ is a cube tiling code and $S=\{a,a',b,b'\}$, then  $m>6\cdot 10^{14}$ (\cite{Ost}). (Let us recall that every cube tilings codes $V,W\subset S^d$ are strongly equivalent for $d\leq 5$ (\cite{DIP,Kii,Ost}), and therefore  $m$ is the number of all cube tiling codes.)}
\end{uwi}

\bigskip
To find covers of a code we shall use the following simple algorithm:

\medskip
{\bf Algorithm} {\small {\sc CoverCode}}.

\medskip
Let $U=\{u^1,...,u^n\}$ be a code, and let $\ka C_{u^i}$, $i\in [n]$, be the family of all covers of the word $u^i\in U$.
Our goal is to find the family $\ka C_U$ of all covers $C_U$ of the code $U$ such that $|C_U|\leq m$ for a fixed $m\in \{1,2...\}$

\medskip   
{\it Input}. The code $U$, the number $m$ and the family $(\ka C_{u^i})_{u^i\in U}$.

{\it Output}. The family $\ka C_U$.

\smallskip
1. For $C_1\in \ka C_{u^1}$ and  $C_2\in \ka C_{u^2}$ if the set $C_1\cup (C_2\setminus C_1)$ is a code (it is obviously a cover of the code $\{u^1,u^2\}$) and has at most $m$ words, then it is attached to the set $\ka C_{1,2}$.  

2. Assuming that the set $\ka C_{1,...,k}$ has already been computed for $2\leq k<n$, we compute the set $\ka C_{1,...,k,k+1}$: for $C \in \ka C_{1,...,k}$ and  $C_{k+1}\in \ka C_{u^{k+1}}$ if the set $C\cup (C_{k+1}\setminus C)$ has at most $m$ words and it is a  code (being a cover of $\{u^1,...,u^{k+1}\}$), then it is attached to the set $\ka C_{1,...,k,k+1}$.

3. $\ka C_U=\ka C_{1,...,n}$. 

\medskip
To present the next algorithm more readable, we define a function  ${\rm Cover}_n(\cdot,\cdot,\cdot)$: 

\smallskip
Let $X,Y\subset S^d$ be polybox codes, $Z\subset S^d$ be a set of words, and let $n\in \en$. (We may think that the position of $Y$ with respect to $X$ is such that $Y$ is not covered by $X$, that is $\bigcup E(Y)\setminus \bigcup E(X)\neq\emptyset$.) 

\smallskip
1. Let $Z_X=\{q\in Z\colon X\cup \{q\}\; {\rm is\;\; a\; code\; and\;\;} |q|_Y>0\}$ (recall that $|q|_Y>0$ means, by \ref{pv}, $\breve{q}\cap \breve{v}\neq\emptyset$ for some $v\in Y$).

\smallskip
2. Let $m=\sum_{v\in Y}(2^d-|v|_X)$ ($m$ can be interpreted as a measure of the uncovered, by $X$, part of $Y$).

\smallskip
3. Let $M$ be the multiset consisting of all numbers $|q|_Y$, $q\in Z_X$, arranged in decreasing order. By $M[i]$ we denote the $i$th element of $M$, and we put $M[i]=0$ if $|M|<i$. Moreover,  $m_n=\sum_{i=1}^nM[i]$.

The number ${\rm Cover}_n(X,Y,Z)\in \{0,1\}$ is defined as follows:
\[ {\rm Cover}_n(X,Y,Z) =
  \begin{cases}
    0       & \quad \text{if } m_n<m,\\
    1  & \quad \text{if } m_n\geq m.
  \end{cases}
\]

\medskip
Let $X\subset S^d$ and $\{u^1,\dots ,u^n\}\subset S^d$ be disjoint polybox codes. If the set of words $U_n(X)=X\cup \{u^1,\dots ,u^n\}$ is a code, then it is called an {\it extension} ({\it by $n$-words}) of $X$. 
Our goal is to find all extensions $U_n(X)$ of $X$ by $n$ words such that each $U_n(X)$ is at the same time a minimal cover of $Y$ (it is assumed that $|v|_Y>0$ for every $v\in X$). Thus, $Z_X$ is the set of all words which can be used to produce such extensions, and ${\rm Cover}_n(X,Y,Z)$ says if there is a chance to do this: If ${\rm Cover}_n(X,Y,Z)=0$, then there is no such chance (the uncovered, by $X$, part of $Y$ is too large); if ${\rm Cover}_n(X,Y,Z)=1$, then it is potentially possible to find such cover. We use ${\rm Cover}_n(X,Y,Z)$ as it can be computed fast and reduces a number configurations to be concerned.

\medskip
{\bf Algorithm} {\small {\sc CoverCode$^\natural$}}.

\medskip   
{\it Input}. Two codes $V$ and $W$ such that $|w|_V>0$ for every $w\in W$ and $W$ is not a cover of $V$; a number $n\in \en$

{\it Output}. The family $\ka C_V^{n+N}$, where $N=|W|$, of all minimal covers $C$ of $V$ such that $W\subset C$, $|C|\leq n+N$ for every $C\in \ka C_V^{n+N}$ and $\delta (C,V)\geq 5$. 

\medskip
\noindent
1. 
Let $Q=\{q\in S^5\colon W\cup \{q\}\; {\rm is\;\; a\; code,\;\;} |q|_V>0\;\;{\rm and}\;\;{\rm Cover}_n(W,V,S^5)=1\}$.

\smallskip
\noindent
2. If $Q=\emptyset$, then $\ka C_V^{n+N}=\emptyset$.

\smallskip
\noindent
3. If $Q\neq\emptyset$, then let 
$$
Q^1=\{W\cup \{q\}\colon q\in Q\;\; {\rm and}\;\; {\rm Cover}_{n-1}(W\cup \{q\},V,Q)=1\}\setminus U^1,\; {\rm where} \; U^1=\{W\cup \{q\}\colon q\in Q \; {\rm and}\;  V\dot{\sqsubseteq} W\cup \{q\}\}.
$$

\smallskip
\noindent
4. For $k\in \{2,...,n-1\}$ we define sets  $Q^k$ and $U^k$: 

\smallskip
If $Q^{k-1}=\emptyset$, then $\ka C_V^{n+N}=U^1\cup\ldots \cup U^{k-1}$.

\smallskip
If $Q^{k-1}\neq\emptyset$, then for every $x\in Q^{k-1}$ let 
$$
Q^k_x=\{x\cup \{q\}\colon q\in Q \;\; {\rm and}\;\; x\cup \{q\} \;\;{\rm is\;\; a\;\; code\;\; and}\;\; {\rm Cover}_{n-k}(x\cup \{q\},V,Q)=1\},
$$
and 
$$
Q^k=\bigcup_{x\in Q^{k-1}}Q^k_x\setminus U^k,\;\;\; {\rm where}\;\;\; U^k=\{x\cup \{q\}\colon x\cup \{q\}\in Q^{k-1} \; {\rm and}\;  V\dot{\sqsubseteq} x\cup \{q\}\}.
$$

\smallskip
\noindent
5. If $Q^{n-1}=\emptyset$, then $\ka C_V^{n+N}=U^1\cup\ldots \cup U^{n-1}$. Otherwise, 
$$
\ka C_V^{n+N}=U^1\cup\ldots \cup U^{n-1}\cup Q^n, \;\; {\rm where}\;\; Q^n=\bigcup_{x\in Q^{n-1}}Q^n_x
$$
and 
$$
Q^n_x=\{x\cup \{q\}\colon q\in Q,\;\; x\cup \{q\} \;\;{\rm is\;\; a\;\; code},\;\;V\sqsubseteq x\cup \{q\}\;{\rm and}\;\delta(x\cup \{q\},V)\geq 5\}\;\; {\rm for}\;\; x\in Q^{n-1}.
$$

\medskip
At some stage of the construction of equivalent codes $V,W$ we know the whole code $W$ and only a portion of $V$, let us denote it by $R$ (that is  $R\sqsubseteq W$). Below, using Lemma \ref{index}, we show how to find the whole $V$. 

\medskip
{\bf Algorithm} {\small {\sc FindSecondCode}}.

\medskip   
{\it Input}. Polybox codes $W,R$ such that $R\sqsubseteq W$, $|R|<|W|$ and $\delta(W,R)\geq 5$.

{\it Output}. The family $\ka V(W)$ of all codes $V$ which are equivalent to $W$, $R\subset V$ and $\Delta(V,W)\geq 5$.

\medskip
1. Let $Q=\{q \in S^d\colon R\cup \{q\} \;\; {\rm is}\;\; {\rm a}\;\; {\rm polybox}\;\; {\rm code },\;\; q\sqsubseteq W,\;\; \delta(W,\{q\})\geq 5\}.$

\smallskip
2. For every $i\in \beta(W)\setminus \beta(R)$ let $V(i)=\{v\in Q\colon \beta(v)=i\}$.

\smallskip
3. If $V(i)=\emptyset$ for some $i\in \beta(W)\setminus \beta(R)$, then $\ka V(W)=\emptyset$.

\smallskip
4. If $V(i)\neq\emptyset$ for every $i\in \beta(W)\setminus \beta(R)$, then let $K=\prod_{i\in \beta(W)}V(i)$, where 
\[ V(i) =
  \begin{cases}
    \{v\in Q\colon \beta(v)=i\}       & \quad \text{if } i\in \beta(W)\setminus \beta(R),\\
    \{r\in R\colon \beta(r)=i\}  & \quad \text{if } i\in \beta(R).
  \end{cases}
\]


\smallskip
5. Let $\ka K=\{\bigcup_{n=1}^{|W|}\{k_n\}\colon k\in K\}$, where $k_n$ is the word standing in $k\in K$ at the $n$th position.

\smallskip
6. $\ka V(W)=\{V\in \ka K\colon V\;\; {\rm is} \;\; {\rm a}\;\; {\rm code}\;\; {\rm and}\;\; \delta(V,W)\geq 5\}$.

\medskip

\subsection{Outline of the proof of Theorem \ref{al} and the results of the computations}

In this subsection we describe the main steps in the proof of Theorem \ref{al}. Moreover, in a series of statements we give the results of basic computations needed in the proof.

Let us note that  if $V\dot{\sqsubseteq} W$ and $W\dot{\sqsubseteq} V$, then, by Corollary \ref{del4}, $\Delta(V,W)\geq 5$. As it was mentioned in Remarks \ref{u2}, a time needed to decide whether $V\dot{\sqsubseteq} W$ may be long. Therefore it is better to check at the beginning (compare {\sc FindSecondCode}) the condition $\Delta(V,W)\geq 5$ (as it is much more faster) and at the very end check the conditions $V\dot{\sqsubseteq} W$ and $W\dot{\sqsubseteq} V$. 

\smallskip
Recall that $P=\{v,w\}$, where $v=bbbbb$ and $w=b'b'b'bb$, and let $T=\{v,u\}$, where $u=b'b'b'b'b'$. By $C_P$ and $C_T$ w denote minimal covers of $P$ and $T$, respectively, and $C_s\subset C_P$ (or $C_s\subset C_T$)  stands for a minimal cover of $s\in \{v,w,u\}$.

Let $\ka C_P$ be the family of all minimal covers of $P$ such that for every $s\in \{v,w\}$ the minimal cover $C_s\subset C_P$ of the word $s$  does not contain twin pairs. For $C_P\in \ka C_P$ let $[C_P]=\{g(C_P)\colon g\in G(S^5)\;{\rm and}\; g(P)=P\}$. Let $\ka N_P$ be the family of all representatives of orbits $[C_P]$, that is $|\ka N_P\cap [C_P]|=1$ for every $C_P\in \ka C_P$.

The first part of the following assumption, which will be made in our proof, steams from Lemma \ref{pbez}.

\medskip
{\bf A}1. $P\subset V$, $C_P\subset W$ and for every $s\in \{v,w\}$ the minimal cover $C_s\subset C_P$ of the word $s$  does not contain twin pairs. Moreover, $C_P\in \ka N_P$.

\medskip
To explain the second part of A1, note that if $C'_P\subset W$, where $C'_P\not\in \ka N_P$ but $C'_P\in [C_P]$, where $C_P\in \ka N_P$, then taking $g\in G(S^5)$ such that $g(C'_P)=C_P$ and $V'=g(V)$, $W'=g(W)$ we obtain $P\subset V'$ (as $g(P)=P$) and $C_P\subset W'$. Thus, $V'$ and $W'$ satisfy A1.

Clearly, by Lemma \ref{ft}, we should consider on more case: $T\subset V$ and $C_T\subset W$. However, as we shall show in Statement \ref{p145}, it is not necessary.

Let $|V|=N$. To simplify computations we shall indicate the family $\ka C_P$ with the assumption \\$|C_v|,|C_w|\leq N-5$ in the case $S=\{a,a',b,b'\}$ and $|C_v|,|C_w|\leq N-2$ in the case $|S|\geq 6$. For the rest cases (that is  $|C_v|>N-5$ or $|C_w|> N-5$ if $S=\{a,a',b,b'\}$ and $|C_v|>N-2$ or $|C_w|> N-2$ if $|S|\geq 6$ for every $C_P\subset W$ such that the minimal covers $C_v,C_w\subset C_P$ of words $v,w$ do not contain a twin pair)  instead of $A1$ we shall assume that

\medskip
{\bf A}2: $v\in V$ and $C_v\subset W$ for $C_v\in \bigcup_{n=N-4}^{N-1}\ka N_v^n$ if $S=\{a,a',b,b'\}$, and $C_v\in \ka N_v^9$ if $|S|\geq 6$,

\medskip
where the family $\ka N_v^n$ consists of all non-isomorphic twin pair free minimal covers $C_v$ of $v=bbbbb$ with $n\in \{N-4,...,N-1\}$ words. (Observe that if $C'_v$ is a minimal cover of $v$ such that $C'_v=g(C_v)$ for some $g\in G(S^5)$, then $g(v)=v$. Thus, in the same way as for $C_P$ we argue that $C_v\in \bigcup_{n=N-4}^{N-1}\ka N_v^n$.)

Let us note that we do not need consider minimal covers $C_v\subset W$ such that $|C_v|=N$. It follows from the following observation:

\begin{lemat}
\label{noN}
Let $V,W\subset S^5$ be two equivalent codes such that $|V|\leq 15$. For every $v\in V$ and $w\in W$ if $C_v\subset W$ and $C_w\subset V$ are twin pair free minimal covers of $v$ and $w$, respectively,  then $|C_v|<|V|$ or $|C_w|<|V|$.
\end{lemat}
\proof
If on the contrary $|C_v|=|C_w|=|V|$, then $V=C_w$ and $W=C_v$, that is $V,W$ are equivalent polybox codes without twin pairs. By Theorem \ref{p12}, there is precisely one, up to isomorphism, such pair of codes, where $|V|\leq 15$. But none of the codes from the family $\ka N_v^{12}$ has the form, up to isomorphism, of one of the codes enumerated in Theorem \ref{p12}.
\hfill{$\square$}

\medskip
As we mentioned at the beginning of Section 3, to find equivalent codes $V$ and $W$ we start with portions of $V$ and $W$ which can be somehow predicted. Denoting these portions by $V_0\subset V$ and $W_0\subset W$, it is seen that our assumptions describe $V_0$ and $W_0$: We have $V_0=P$ and $W_0=C_P$ or $V_0=\{v\}$ and $W_0=C_v$. The only reason for considering two systems of initial conditions $V_0,W_0$ (that is A1 or A2), and not only one, is a simplification of our computations: If we resign of the assumption A2, and in A1 we assume that $|C_v|,|C_w|\leq N-1$, then the family $\ka C_P$ is much bigger than in the case $|C_v|,|C_w|\leq N-5$, and indication of $\ka C_P$ is much more longer. On the other hand, if we consider only initial conditions $v\in V$ and $C_v\subset W$, where $|C_v|\in \{5,...,14\}$, then reconstruction of $W$ in the base of $C_v$ for small numbers $|C_v|$ is difficult (for example, if $|C_v|=5$ and $|W|=15$, then we have to extend $C_v$ by 10 words to get $W$).

\smallskip
Our programs were written in Julia and Python. The total time for computations was around three days (we used 8-core 3.4-GHz processor with 32GB RAM memory).

\smallskip
To find the family $\ka N_P$ and $\bigcup_{n=N-4}^{N-1}\ka N_v^n$, we have to compute twin pairs free covers of a word. Using algorithm {\sc CoverWord} we obtained the following results:

\begin{st}
\label{sl11}
Let $S=\{a,a',b,b'\}$, $v=bbbbb$, and let $\ka N^n_v$ be a family of all non-isomorphic twin pair free minimal covers of $v$ such that $|C_v|\in \{5,...,14\}.$ 
Then
$$
|\ka N^5_v|=1,\;\;|\ka N^6_v|=1,\;\; |\ka N^7_v|=3,\;\; |\ka N^8_v|=4,\;\; |\ka N^9_v|=19,\;\; |\ka N^{10}_v|=51,
$$
$$
|\ka N^{11}_v|=153,\;\;|\ka N^{12}_v|=287,\;\;|\ka N^{13}_v|=683,\;\;|\ka N^{14}_v|=1275.
$$

Additionally, if $S=\{a,a',b,b',c,c'\}$, then there are four twin pair free minimal covers $C_v$ of $v$ with $|C_v|=9$ such that $C_v=C_v^{i,a}\cup C_v^{i,a'}\cup C_v^{i,b}\cup C_v^{i,c}\cup C_v^{i,c'}$ for some $i\in [5]$, where all sets $C_v^{i,s}$, $s\in \{a,a',c,c'\}$, are non-empty.
\hfill{$\square$}
\end{st}

\smallskip
The following result of computations shows that in A1 we do not need consider the case $T\subset V$:

\begin{st}
\label{p145}
Let $S=\{a,a',b,b'\}$, $T=\{v,u\}$, where $v=bbbbb$ and $u=b'b'b'b'b'$. Let $\ka C^0_s$, $s\in \{v,u\}$, be the family of all minimal twin pair free covers $C^0_s\subset S^5$ of $s\in \{v,u\}$ such that $|C^0_s|\in \{5,...,14\}$ and neither $C^0_v$ nor $C^0_u$ contains a pair of words which is isomorphic to the pair $\{bbbbb,b'b'b'bb\}$. If $\ka C_T$ is the family of all minimal covers $C_T\subset S^5$ of $T$ computed according algorithm {\sc CoverCode}, in which we take $\ka C_{u^1}=\ka C^0_v$ and $\ka C_{u^2}=\ka C^0_u$, such that $|C_T|\leq 15$ and $C_T$ does not contain a pair of words which is isomorphic to the pair $\{bbbbb,b'b'b'bb\}$ for every $C_T\in \ka C_T$, then $\ka C_T=\emptyset$. Moreover, the same is true if $|S|\geq 6$, $|C^0_s|\in \{5,...,9\}$ for $s\in \{v,u\}$ and $|C_T|\leq 10$ for $C_T\in \ka C_T$.
\hfill{$\square$}
\end{st}

In the next two statements, we give families of covers of the code $P$. In the computations we used results presented in Statement \ref{sl11} and algorithm {\sc CoverCode}.

\begin{st}
\label{sabc}
Let $S=\{a,a',b,b',c,c'\}$, $P=\{v,w\}$, where $v=bbbbb$ and $w=b'b'b'bb$. Let $\ka C_s$, $s\in \{v,w\}$, be the family of all minimal twin pair free covers $C_s\subset S^5$ of $s\in \{v,w\}$ such that $|C_s|\in \{5,6,7,8\}$, and let $\ka C_P$ be the family of all minimal covers $C_P\subset S^5$ of $P$ computed according algorithm {\sc CoverCode}, in which we take $\ka C_{u^1}=\ka C_v$ and $\ka C_{u^2}=\ka C_w$, such that $|C_P|\leq 10$ for $C_P\in \ka C_P$.
Then $|\ka C_P|=214080$ and 
$$
|\ka C^8_P|=64,\;\;\; |\ka C^9_P|=1536,\;\;\;|\ka C^{10}_P|=212480.
$$
Moreover, 
$$
|\ka N^8_P|=0,\;\;\; |\ka N^9_P|=0,\;\;\;|\ka N^{10}_P|=58,
$$
where $\ka N^m_P$, $m \in \{8,9,10\}$, consists of all codes $C_P\in \ka N_P$ such that $|C_P|=m$  and $P\dot{\sqsubseteq} C_P$.
\hfill{$\square$}
\end{st}

\begin{st}
\label{p14}
Let $S=\{a,a',b,b'\}$, $P=\{v,w\}$, where $v=bbbbb$ and $w=b'b'b'bb$. Let $\ka C_s$, $s\in \{v,w\}$, be the family of all minimal twin pair free covers $C_s\subset S^5$ of $s\in \{v,w\}$ such that $|C_s|\in \{5,...,10\}$, and let $\ka C_P$ be the family of all minimal covers $C_P\subset S^5$ of $P$ computed according algorithm {\sc CoverCode}, in which we take $\ka C_{u^1}=\ka C_v$ and $\ka C_{u^2}=\ka C_w$, such that $|C_P|\leq 15$ for $C_P\in \ka C_P$.
Then $|\ka C_P|=4965112$ and 
$$
|\ka C^8_P|=8,\;\;\; |\ka C^9_P|=96,\;\;\;|\ka C^{10}_P|=4256,\;\;\;|\ka C^{11}_P|=17760,\;\;\;|\ka C^{12}_P|=158048,\;\;\;|\ka C^{13}_P|=449568,
$$
$$
|\ka C^{14}_P|=1795552,\;\;\;|\ka C^{15}_P|=2539824.
$$
Moreover, 
$$
|\ka N^8_P|=0,\;\;\; |\ka N^9_P|=0,\;\;\;|\ka N^{10}_P|=23,\;\;\;|\ka N^{11}_P|=42,\;\;\;|\ka N^{12}_P|=379,\;\;\;|\ka N^{13}_P|=839,
$$
$$
|\ka N^{14}_P|=3679,\;\;\;|\ka N^{15}_P|=3665,
$$
where $\ka N^m_P$, $m \in \{8,...,15\}$, consists of all codes $C_P\in \ka N_P$ such that $|C_P|=m$ and $P\dot{\sqsubseteq} C_P$. 
\hfill{$\square$}
\end{st}

At the end, we give a result which simplify the proof of the first part of Theorem \ref{al} in the case $|S|\geq 8$.

\begin{lemat}
\label{cd}
Let $v=bbbbb$, $w=b'b'b'bb$, $P=\{v,w\}$, and let $C_P\subset S^5$ be a minimal cover of $P$ such that minimal covers $C_v,C_w\subset C_P$ of $v$ and $w$, respectively do not contain a twin pair. If $|C_P|\leq 8$, then the code $C_P$ can be written down in the alphabet $S=\{a,a',b,b',c,c'\}$.
\end{lemat}
\proof
Since $|C_P|\leq 8$, we have $|C_s|\leq 8$ for $s\in \{v,w\}$. There are only 9 non-isomorphic  twin pairs free covers of the word $v$ with $|C_v|\leq 8$ and for each of them we have $C_v\subset \{a,a',b\}^5$. Similarly  $C_w\subset \{a,a',b,b'\}^5$. Thus, $C_P=C_P^{i,a}\cup C_P^{i,a'}\cup C_P^{i,b}\cup C_P^{i,b'} \cup C_P^{i,c}\cup C_P^{i,c'}$ for $i\in [5]$. Indeed, if $C_v\subset \{a,a',b\}^5$, then $C_w\subset \{a,a',b,b'\}^5$ or $C_w\subset \{c,c',b,b'\}^5$. Clearly, $C_P=C_v\cup C_w$.
\hfill{$\square$}

\section{Proofs}

We divide the proof of Theorem \ref{al} with respect to the number of letters in $S$.

\medskip
Let us recall that if  $C\subset S^d$ and $\{u^1,\dots ,u^n\}\subset S^d$ are disjoint polybox codes, and $U_n(C)=C\cup \{u^1,\dots ,u^n\}$ is a code, then $U_n(C)$ is called an extension (by $n$-words) of $C$.  Let $\ka U_n(C)$ be the family of all extensions $U_n(C)$ of the code $C$. 

\subsection{The case $|V|\leq 10$ and $|S|\geq 6$}

\medskip

{\it Proof of the first part  of Theorem \ref{al}}. Let $V,W\subset S^5$ be equivalent codes which satisfy A1 or A2 (compare the previous section). Since $V\dot{\sqsubseteq} W$ and $W\dot{\sqsubseteq} V$ we have, by Corollary \ref{del4}, $\Delta(V,W)\geq 5$.  Let $N=|V|$, and let $\ka N^m_P,\ka N^m_v$ be as in  Statements \ref{sabc} and \ref{sl11}, respectively. 

If $|S|\geq 8$ and $|V|\leq 8$, then, by Lemma \ref{cd}, we may assume that $V\subset \{a,a',b,b',c,c'\}^5$.

By Statement \ref{sabc}, there are no covers in $\ka C_P$ with less that 8 words. This proves our theorem in the case $N\leq 7$ (that is, every equivalent codes $V,W\subset S^5$ with $|V|\leq 7$ are strongly equivalent).
 
Since, again by Statement \ref{sabc}, $\ka N^8_P=\emptyset$ the theorem, by Lemma \ref{ind}, is also true in the case $N=8$. 

Let $N=9$. By Lemma \ref{noN}, we do not need to consider covers from $\ka N^9_v$, and since $\ka N^9_P=\emptyset$, the theorem holds true in the case $N=9$.


Let $N=10$, and let $\ka U^9_1=\bigcup_{C\in \ka N^9_v}\ka U_1(C)$.

We computed that
$$
|\ka U^9_1|=18382.
$$
For every $W\in \ka U^9_1$, using algorithm {\sc FindSecondCode} in which $R=\{v\}$, we obtained $\ka V(W)=\emptyset$. 

Similarly, $\ka V(W)=\emptyset$ for every $W\in \ka N^{10}_P$, where in algorithm {\sc FindSecondCode} we take $R=P$. Since $\ka N^8_P\cup \ka N^9_P=\emptyset$, by Lemma \ref{ind}, the theorem is proved.
\hfill{$\square$}

\medskip
Since in the above proof we considered also flat codes (which can be identify with codes in dimensions $d<5$), we have the following corollary:

\begin{wn}
\label{md}
Every equivalent codes $V,W\subset S^d$, where $d\leq 5$ and $|V|\leq 8$ are strongly equivalent. \hfill{$\square$}
\end{wn}

\subsection{The case $|V|\leq 15$ and $S=\{a,a',b,b'\}$}

Some computations made in our proof of the second part of Theorem \ref{al} deal with flat polybox codes. For example, we shall compute extensions by $n\in \{0,...,4\}$ words of codes $W$ with $|W|=11$ and $W=W^{i,a}$ for some $i\in [5]$. For some $n$ (especially for $n=4$) such computations are long (because of the flatness of a code), and since we are able to predict the properties of the resulting extensions, some of the computations can be omitted, and some can be simplify. Below we give results which help us to do such simplifications. 

We start with an interesting observation steaming from the first part of Theorem \ref{al}. 

Let $V\subset S^d$ be a cube tiling code, $Q\subset V$ be a simple component of $V$, and let $U(Q)$ be the simple cube tiling code containing $Q$. If it is possible to pass from $V$ to $U(Q)$ by gluing and cutting such that the code $Q$ is unchanged during the process of gluing and cutting, then $Q$ is called {\it fixed} simple component.

\begin{tw} 
\label{stab4}
Let $V\subset S^4$ be a cube tiling code, $Q\subset V$ be a simple component of $V$, and let $U(Q)$ be the simple cube tiling code containing $Q$. Then it is possible to pass from $V$ to $U(Q)$ by gluing and cutting, where $Q$ is fixed simple component.  
\end{tw}
\proof
Let $S=\{a_1,a_1',...,a_8,a_8'\}$. For simplicity we may assume that $U(Q)=\{a_1,a_1'\}^4$. Let $i\in [4]$ be such that $V=V^{i,a_1}\cup V^{i,a_1'}\cup \ldots \cup V^{i,a_n}\cup V^{i,a_n'}$ for some $n\in \{2,...,8\}$, where all sets on the right side of the decomposition of $V$ are non-empty. (If $V=V^{i,a_1}\cup V^{i,a_1'}$ for every $i\in [4]$, then $V=U(Q)$ and there is nothing to prove.) Since, by Corollary \ref{md}, for every $j\in \{2,...,n\}$, the codes $V_{i^c}^{i,a_j}$ and $V_{i^c}^{i,a'_j}$ are strongly equivalent, we may pass by gluing and cutting from $V$ to a cube tiling code $\bar{V}$, where $\bar{V}=\bar{V}^{i,a_1}\cup \bar{V}^{i,a_1'}$ (compare Subsection 2.2). Note that, during the process of passing from $V$ to $\bar{V}$ the code $V^{i,a_1}\cup V^{i,a_1'}$ was unchanged, and thus $Q$ was unchanged, as $Q\subset  V^{i,a_1}\cup V^{i,a_1'}$. In this manner we can pass from $V$ to $U(Q)$ by gluing and cutting keeping $Q$ unchanged.
\hfill{$\square$}

\begin{wn} 
\label{w15}
Every equivalent polybox codes $V,W\subset S^4$ with $|V|=15$ are strongly equivalent.
\end{wn}
\proof
In \cite{DIP} it was shown that $V$ and $W$ can be extended to a cube tiling code by one word. Denote this word by $v\in S^4$. Let $Q_V\subset V\cup \{v\}$ and $Q_W\subset W\cup \{v\}$ be the simple components of $V\cup \{v\}$ and $W\cup \{v\}$, respectively such that $v\in Q_V\cap Q_W$. Let $U(Q_V)$ be the simple cube tiling code containing $Q_V$. By Theorem \ref{stab4}, we can pass from $V\cup \{v\}$ and $W\cup \{v\}$ to $U(Q_V)$ by gluing and cutting with $Q_V$ and $Q_W$ fixed. In particular, the word $v$ is unchanged during the process of gluing and cutting. Thus, $V$ and $W$ are strongly equivalent.
\hfill{$\square$}

\medskip

{\it Proof of the second part of Theorem \ref{al}.} Let $V,W\subset S^5$ be equivalent codes satisfying A1 and A2 (see Subsection 3.2). As in the proof of the first part, we have $\Delta(V,W)\geq 5$. Let $N=|V|$, and let $\ka N^m_P,\ka N^m_v$ be as in  Statements \ref{p14} and \ref{sl11}, respectively. By the first part of Theorem \ref{al} every equivalent codes $V,W\subset S^5$ with $|V|\leq 10$ are strongly equivalent and thus, we have to consider the case $11\leq N\leq 15$. 

Let $\ka N=\ka N^{10}_P\cup \ka N_P^{15}\cup \bigcup_{i=11}^{14}\ka N_v^{i}\cup \ka N_P^{i}$. In the first step we shall compute all extensions $\ka U_n(C)$ of codes $C\in \ka N$ for $n\in \{0,...,4\}$, where $|W|\in \{11,...,15\}$ for every $W\in \ka U_n(C)$. Note that, if some $C$ is flat, that is $C=C^{i,s}$ for some $i\in [5]$ and $s\in S$ then, by Lemma \ref{fla}, computing $U_n(C)$ we can extend $C$ by words $v$ with $v_i=s$, which speeds up computations. (For $C\in \ka N_P^{10}$  extensions $U_5(C)$ will be computed in a little different manner than the others extensions.) Next, using algorithm {\sc FindSecondCode}, for every $W\in \bigcup_{n=0}^4\bigcup_{C\in \ka N}\ka U_n(C)$ we compute the family $\ka V(W)$ of codes $V$ which are  equivalent to $W$ and $\Delta(V,W)\geq 5$. In Table 1 we collected the numbers of such pairs 
with one exception: We did not compute extensions $\ka U_4(C)$ for four codes $C\in \ka N_P^{11}$. The reason is that all these four codes are flat, and thus they can be identify with codes in dimension four. Their extensions $W\in \ka U_4(C)$ have 15 words. Since we are looking for equivalent codes $V,W$ with $\Delta(V,W)\geq 5$, by Lemma \ref{fla}, extensions $W\in \ka U_4(C)$ have to by flat. Thus, by Corollary \ref{w15}, for every $W\in \ka U_4(C)$ and every code $V$ which is equivalent to $W$ the codes $V,W$ are strongly equivalent. (In Table 1 we denoted the lack of these computations by $0^\dag$, that is for mentioned four codes we did not make any computations, and for the rest (non-flat) codes $W\in \ka U_4(C)$ we found no pairs $V,W$ of equivalent codes such that $|V|=15$ and $\Delta(V,W)\geq 5$.)

We should also find extensions of $C_P\in \ka N_P^{10}$ by five words. However, such extensions need long computations (we have to add five words), and therefore we shall use algorithm {\sc CoverCode$^\natural$}, which is much more faster in this case.  To do this, recall that we assume that $P\subset V$ and $C_P\subset W$, where covers $C_v,C_w\subset C_P$ of $v$ and $w$ do not contain a twin pair. 
An inspection of the set $\ka N_P^{10}$ shows that each $C_P\in \ka N_P^{10}$ contains a pair of words $Q=\{p,q\}$ which is isomorphic to $P$. Let $C_Q\subset V$ be  a minimal cover of $Q$. By Lemma \ref{bez} we pass  by gluing and cutting to a code $\bar{C}_Q$  in which covers $C_p,C_q\subset \bar{C}_Q$ do not contain a twin pair and $|\bar{C}_Q|\leq 10$ (we computed at the beginning extensions of codes from $\ka N_P^{i}$ for $i \in \{11,...,14\}$). Additionally, we may exclude from further consideration all these covers $\bar{C}_Q$ for which one can pass by gluing and cutting to a code $\bar{\bar{C}}_Q$  such that $\bar{\bar{C}}_Q\cap Q\neq\emptyset.$ Clearly, passing from $C_Q$ to $\bar{C}_Q$ changes $V$ to $\bar{V}$. Thus, we obtain a pair of equivalent codes $\bar{V}$ and $W$ such that $C_P\subset W$, with $C_P\in \ka N^{10}_P$ and $\bar{C}_Q\subset \bar{V}$, where $\bar{C}_Q\in \ka C^{10}_Q$ and $\ka C^{10}_Q$ is the family of all covers of $Q$ with 10 words (let us recall that, by Statement \ref{p14}, $\ka N_Q^i=\emptyset$ for $i\in [9]$). Now, using algorithm {\sc CoverCode$^\natural$}, we compute all extensions $U^\natural_5(C_P)$ of $C_P$ with the assumption that these extensions have to cover a code $\bar{C}_Q$. Denote the family of all $U^\natural_5(C_P)$ by $\ka U^\natural_5(C_P)$. Now for every  $C_P\in \ka N_P^{10}$ and every $W\in \ka U^\natural_5(C_P)$, using algorithm {\sc FindSecondCode} in which we take $R=\bar{C}_Q$, we compute the family $\ka V(W)$. We obtained 8 pairs $W,\bar{V}$, where $\bar{V}\in \ka V(W)$.

\smallskip
\begin{center}
\begin{tabular}{c| c| c | c | c |c }
\hline
    &  $|U(C)|=15$ &   $|U(C)|=14$   &  $|U(C)|=13$ &   $|U(C)|=12$  &  $|U(C)|=11$ \\

\hline
$C\in \ka C_P^{10}$ &  $8^\dag$ & 46 &  0 & 0& 0  \\
$C\in\ka N_v^{11}$ &  1 & 0 &  0 & 0& 0 \\
$C\in\ka N_P^{11}$ &  $0^\dag$ & 20 &  1 & $2^\dag$& 0 \\
$C\in\ka N_v^{12}$ &  0 & 0 &  0 & 0& - \\
$C\in\ka N_P^{12}$ &  0 & 0 &  0 & 0& - \\
$C\in\ka N_v^{13}$ &  0 & 0 &  0 & -& - \\
$C\in\ka N_P^{13}$ &  0 & 0 &  0 & -& - \\
$C\in\ka N_v^{14}$ &  0 & 0 &  - & -& - \\
$C\in\ka N_P^{14}$ &  0 & 0 &  - & -& - \\
$C\in\ka N_P^{15}$ &  0 & - &  - & -& - \\

\end{tabular}
\end{center}

\smallskip
\noindent{\footnotesize Table 1: The number of pairs of equivalent codes $V,W$ ($|W|=|U(C)|$, $U(C)$ is an extension of $C$ by $n\in \{0,...,4\}$ words and $V\in \ka V(W)$) with $\Delta(V,W)\geq 5$ and $11\leq |V|\leq 15$. $8^\dag$ means that there are 8 pairs $\bar{V},W$ found using algorithm {\sc CoverCode$^\natural$}; $0^\dag$ means that we did not make computations for flat codes $C\in \ka N_P^{11}$, and for the rest (non-flat) codes we found no pairs with desire properties. Finally, $2^\dag$ means that one of the two pairs is the special pair.}

\smallskip

Now for every computed pairs of codes $V,W$ and $\bar{V},W$ we have to check whether $V\dot{\sqsubseteq} W$, $W\dot{\sqsubseteq} V$ and $\bar{V}\dot{\sqsubseteq} W$, $W\dot{\sqsubseteq} \bar{V}$, respectively. To do this, we check first whether it is possible to pass by gluing and cutting from $W$ to $\bar{W}$ such that $\bar{W}\cap P\neq\emptyset$ (in the case of codes $W,\bar{V}$ we check whether it is possible to pass by gluing and cutting from $\bar{V}$ to $\bar{\bar{V}}$  such that $\bar{\bar{V}}\cap Q\neq\emptyset$, where $Q\subset W$).

We computed, using algorithm {\sc GlueAndCut}, that for each considered pair of codes, except one, such passing is possible, that is it is not true that $V\dot{\sqsubseteq} W$, $W\dot{\sqsubseteq} V$ and $\bar{V}\dot{\sqsubseteq} W$, $W\dot{\sqsubseteq} \bar{V}$ for all, but one, pairs $V,W$ and $\bar{V},W$. This exceptional pair $V,W$ is, up to isomorphism, the special pair. In particular, by the first part of the theorem and Lemma \ref{ind}, we showed that every equivalent codes $V,W\subset S^5$, $S=\{a,a',b,b'\}$, with $|V|\leq 11$ are strongly equivalent. 

\hfill{$\square$}

\begin{wn}
\label{dec}
Let $V,W\subset S^d$, $S=\{a,a',b,b'\}$, $d\leq 5$, be two equivalent polybox codes which are not strongly equivalent. Then $d\geq 4$, and for $d=4$ the codes $V,W$ form, up to isomorphism, the special pair. If $d=5$ and $|V|\leq 15$, then $|V|\in \{12,...,15\}$ and it is possible to pass from $V$ to $\bar{V}$  and from $W$ to $\bar{W}$ by gluing and cutting such that there are  decompositions $\bar{V}=V^1\cup P$ and $\bar{W}=W^1\cup P$, where $V^1,W^1$ is, up to isomorphism, the special pair.
\end{wn}
\proof
By Lemma \ref{sind} and Theorem \ref{al} we have to prove only the first part of the corollary. Since, by Theorem \ref{al}, the special pair $V,W$ is the only pair of equivalent codes $V,W\subset S^d$, where $d\leq 5$ and $|V|\leq 15$, with $V\dot{\sqsubseteq} W$ and  $W\dot{\sqsubseteq} V$, we have $d\geq 4$. It is easy to check that there are only four words $V'=\{v^1,...,v^4\}$ such that the sets $V\cup Q$, $W\cup Q$ are polybox codes for any $Q\subset V'$, where $V,W\subset S^4$ is the special pair. Moreover, it can be computed (using algorithm {\sc GlueAndCut}) that for every non-empty set $Q\subset V'$ the codes $V\cup Q$ and $W\cup Q$ are strongly equivalent. Thus, by Lemma \ref{sind}, the special pair is the only pair, up to isomorphism, of equivalent codes in dimension four which are not strongly equivalent.
\hfill{$\square$}

\subsection{Proof of Theorem \ref{gc}}

{\it Proof of Theorem \ref{gc}}. Let $S=\{a_1,a_1',...,a_k,a_k'\}$, and let $i\in [6]$ be such that the number of sets on the right side of the decomposition
$$
U=U^{i,a_1}\cup U^{i,a_1'}\cup \ldots \cup U^{i,a_n}\cup U^{i,a_n'}
$$
is the greatest, $n\leq k$, and  all these sets are non-empty. If $n\geq 4$, then we may assume that $|U^{i,a_n}\cup U^{i,a_n'}|\leq 16$. Let $V=U^{i,a_n}_{i^c}$ and $W=U^{i,a'_n}_{i^c}$. By the first part of Theorem \ref{al}, the codes $V,W\subset S^5$ are strongly equivalent. Thus, proceeding as described in Subsection 2.2 we can pass from $U$ to a cube tiling code $\bar{U}$ by gluing and cutting, where 
$$
\bar{U}=U^{i,a_1}\cup U^{i,a_1'}\cup \ldots \cup U^{i,a_{n-2}}\cup U^{i,a_{n-2}'} \cup \bar{U}^{i,a_{n-1}}\cup \bar{U}^{i,a_{n-1}'}.
$$
In this way we can reduce the number of sets in the decomposition of words in $U$, and therefore we can assume, again by the first part of Theorem \ref{al}, that it is possible to pass from $U$ to $\bar{U}\subset \{a,a',b,b'\}^6$ by gluing and cutting (we took $a_1=a,a_2=b$). 

Suppose that
$$
\bar{U}=\bar{U}^{i,a}\cup \bar{U}^{i,a'}
$$
for some $i\in [6]$, and let $V=\bar{U}^{i,a}_{i^c}$ and $W=\bar{U}^{i,a'}_{i^c}$. Since $|V|=|W|=32$, the codes $V,W\subset S^5$ are two cube tiling codes. Let $P\subset S^5$ be any fixed simple cube tiling code.  As every two cube tilings codes in dimension five are strongly equivalent  (\cite{Ost}), it follows that we may pass from $V$ to $P$ and similarly, from $W$ to $P$ by gluing and cutting. Thus, by gluing and cutting we can pass from $\bar{U}$ to $\bar{\bar{U}}$, where 
$$
\bar{\bar{U}}=\bar{\bar{U}}^{i,a}\cup \bar{\bar{U}}^{i,a'},
$$
and $\bar{\bar{U}}^{i,a}_{i^c}=P$, $\bar{\bar{U}}^{i,a'}_{i^c}=P$, that is $\bar{\bar{U}}$ is a simple cube tiling code.

Let us consider now the case in which 
\begin{equation}
\label{y}
\bar{U}=\bar{U}^{i,a}\cup \bar{U}^{i,a'}\cup \bar{U}^{i,b}\cup \bar{U}^{i,b'}
\end{equation}
for every $i\in [6]$. We shall show that we may assume that $|\bar{U}^{i,s}\cup \bar{U}^{i,s'}|\leq 30$ for some $i\in [6]$ and $s\in \{a,b\}$.

By \cite{P}, the cube tiling code $\bar{U}$ contains a twin pair. Without loss of generality, we may assume that $v$ and $w$ form this pair, where $v_i=b$ and $w_i=b'$, that is $v\in \bar{U}^{i,b}$ and  $w\in \bar{U}^{i,b'}$. If $|\bar{U}^{i,b}\cup \bar{U}^{i,b'}|>32$, then $|\bar{U}^{i,a}\cup \bar{U}^{i,a'}|\leq 30$ (the number $|\bar{U}^{i,s}\cup \bar{U}^{i,s'}|$ is even) . Let $|\bar{U}^{i,b}\cup \bar{U}^{i,b'}|=32$. Passing from the code $\{v,w\}$ by gluing and cutting to the twin pair $\{\bar{v},\bar{w}\}$, where $\bar{v}_i=a$ and $\bar{w}_i=a'$, we pass from $\bar{U}$ to $\bar{\bar{U}}$ by gluing and cutting, where $|\bar{\bar{U}}^{i,b}\cup \bar{\bar{U}}^{i,b'}|=30$.Thus, in (\ref{y}) we may assume that $|\bar{U}^{i,b}\cup \bar{U}^{i,b'}|\leq 30$.


Let $V=\bar{U}^{i,b}_{i^c}$ and $W=\bar{U}^{i,b'}_{i^c}$. Clearly, the codes $V,W\subset S^5$ are equivalent and $|V|\leq 15$. 

If $V$ and $W$ are strongly equivalent, then we can pass (compare Subsection 2.2) from $\bar{U}$ to $\bar{\bar{U}}$ by gluing and cutting, where $\bar{\bar{U}}$ is such that 
$$
\bar{\bar{U}}=\bar{\bar{U}}^{i,a}\cup \bar{\bar{U}}^{i,a'}.
$$
This means, as we showed above, that we can pass by gluing and cutting from $U$ to a simple cube tiling code.

Assume now that $V$ and $W$ are not strongly equivalent. By Corollary \ref{dec}, it is possible to pass from $V$ to $\bar{V}$ and from $W$ to $\bar{W}$ by gluing and cutting, where $\bar{V}=\bar{V}^1\cup P$, $\bar{W}=\bar{W}^1\cup P$ and the pair of codes $\bar{V}^1,\bar{W}^1$ is, up to isomorphism, the special pair in $S^5$. We may assume that $i=6$ and $v_5=a$ for every $v\in \bar{V}^1\cup \bar{W}^1$ (compare Theorem \ref{p12}). Since $|\bar{V}^1\cup \bar{W}^1|=24$, it follows that $|\bar{\bar{U}}^{5,a}|\geq 24$. 
Then $|\bar{\bar{U}}^{5,b}\cup \bar{\bar{U}}^{5,b'}|\leq 16$ and consequently, by Theorem \ref{al}, the codes $\bar{\bar{U}}^{5,b}_{5^c}$ and $\bar{\bar{U}}^{5,b'}_{5^c}$, are strongly equivalent. Therefore, in the same manner as above we show that we can pass by gluing and cutting from $U$ to a simple cube tiling code.
\hfill{$\square$}

\medskip
\begin{uw}{\rm
The above proof shows that using Theorem \ref{al} we can prove Theorem \ref{gc} independently from the result given in \cite{Ost}.}
\end{uw}

\section{Some consequences and open problems}

In the last section we give a few results that steam mainly from Theorem \ref{al}. We believe that there are much more interesting problems dealing with glue and cut procedure. At the end of the paper we shall formulate some of them. 

\subsection{Some consequences of Theorem \ref{al}}

Theorem \ref{stab4} can be extended on cube tiling codes in dimension five:

\begin{tw} 
\label{stab5}
Let $V\subset S^5$ be a cube tiling code, $Q\subset V$ be a simple component of $V$, and let $U(Q)$ be the simple cube tiling code containing $Q$. Then it is possible to pass from $V$ to $U(Q)$ by gluing and cutting with  $Q$ fixed.  
\end{tw}
\proof
Suppose that there is $i\in [5]$ such that for every $j\in [k]$ codes $V^{i,a_j}_{i^c}, V^{i,a_j'}_{i^c}$ are strongly equivalent, where $V=V^{i,a_1}\cup V^{i,a_1'}\cup \cdots \cup V^{i,a_{k}}\cup V^{i,a_{k}'}$. We may assume that $Q\subset V^{i,a_1}\cup V^{i,a'_1}$.  Since the codes $V^{i,a_j}_{i^c}$ and $V^{i,a'_j}_{i^c}$ are strongly equivalent for $j\in \{2,...,k\}$ we can pass by gluing and cutting  from $V$ to $\bar{V}$, where $\bar{V}=\bar{V}^{i,a_1}\cup \bar{V}^{i,a'_1}$ and $ V^{i,a_1}\cup V^{i,a'_1}\subset \bar{V}$, which means that $Q\subset \bar{V}$. Let $Q^1=Q^{i,a_1}_{i^c}, Q^2=Q^{i,a'_1}_{i^c},V^1=\bar{V}^{i,a_1}_{i^c}$ and $V^2=\bar{V}^{i,a'_1}_{i^c}$. By  Theorem \ref{stab4}, it is possible pass from $V^1$ to $U_{i^c}(Q)$ and from $V^2$ to $U_{i^c}(Q)$ by gluing and cutting such that $Q^1$ and $Q^2$ are fixed. Hence, we can pass from $V$ to $U(Q)$ by gluing and cutting with $Q$ fixed.

Now assume that there is $r\in [5]$, we may assume that $r=5$, such that  the codes $V^{5,a_1}_{5^c}, V^{5,a_1'}_{5^c}$ are not strongly equivalent. Then, by Corollary \ref{dec} (for $d=4$), the codes $V^{5,a_1}_{5^c}, V^{5,a_1'}_{5^c}$ form , up to isomorphism, the special pair. As it was mentioned in the proof of Corollary \ref{dec}, there are four words $v^1,...,v^4\in S^4$ such that  $V^{5,a_1}_{5^c}\cup \{v^1,...,v^4\}$ and $V^{5,a'_1}_{5^c}\cup \{v^1,...,v^4\}$ are cube tiling codes. Moreover, no two words from the set $\{v^1,...,v^4\}$ belong to the same simple code and the codes $V^{5,a_1}_{5^c}\cup P$ and $V^{5,a'_1}_{5^c}\cup P$ are strongly equivalent for every nonempty $P\subset \{v^1,...,v^4\}$. It follows that 
$V=V^{5,a_1}\cup V^{5,a_1'}\cup \{v^1b_1,v^1b'_1\}\cup \ldots \cup \{v^4b_4,v^4b'_4\}$,  where $b_1,...,b_4\in S$, 
and we may assume that $Q\cap \{v^1b_1,v^1b'_1\}=\emptyset$. Passing from $\{v^1b_1,v^1b'_1\}$ to $\{v^1a_1,v^1a'_1\}$ by gluing and cutting, we pass from $V$ to $\bar{V}$, where $\bar{V}$ is, by the first part of Corollary \ref{dec}, such as at the beginning of the proof.
\hfill{$\square$}

\medskip
Theorem \ref{stab5} is a source of non-trivial examples of strongly equivalent codes $V,W\subset S^5$:

\begin{wn} 
\label{eq5}
Let $V,W\subset S^5$ be equivalent polybox codes and suppose that there is a simple code $V'\subset S^5$ such that $V\cup V'$ is a cube tiling code. Then $V$ and $W$ are strongly equivalent. In particular,  if $V$ is a cube tiling code, $Q\subset V$ is a simple component of $V$ and $W$ is equivalent to $V\setminus Q$, then $V\setminus Q$ and $W$ are strongly equivalent.
\hfill{$\square$}
\end{wn}

From Corollary \ref{eq5} we obtain immediately the following two observations:

\begin{wn} 
\label{nst5}
Let $V,W\subset S^5$ be equivalent polybox code which are not strongly equivalent. Then for any simple codes $V',W'\subset S^5$, the sets $V\cup V'$ and $W\cup W'$ are not  cube tiling codes.  
\hfill{$\square$}
\end{wn}

\begin{wn} 
\label{nstt5}
If $V\subset S^5$ is a cube tiling code which is not simple, and $Q\subset V$ is a simple component of $V$, then the polybox code $V\setminus Q$ contains a twin pair.   
\end{wn}
\proof
Let $U(Q)$ be the simple cube tiling code containing $Q$. The codes $V\setminus Q$ and $U(Q)\setminus Q$ are equivalent. Note, that  $V\setminus Q\cap U(Q)\setminus Q=\emptyset$, otherwise $Q$ would not be a simple component of $V$. By Corollary \ref{eq5}, the codes $V\setminus Q$ and $U(Q)\setminus Q$ are strongly equivalent, and since $V\setminus Q\neq U(Q)\setminus Q$, the code $V\setminus Q$ contains a twin pair.
\hfill{$\square$}

\medskip
A polybox code $V\subset S^d$ is {\it rigid} if for every code $W\subset S^d$ which is equivalent to $V$ we have $V=W$ (\cite{KP}).

\begin{wn}
\label{szt}
Every polybox code $V\subset S^d$ without twin pairs, where $S=\{a,a',b,b'\}$ and $d\leq 5$, having at most 11 words is rigid. 
\end{wn}
\proof
Assume on the contrary that $V$ is not rigid, and let $W\neq V$ be a code which is equivalent to $V$. By Theorem \ref{al}, $V$ and $W$ are strongly equivalent, and thus both codes contain twin pairs. A contradiction.
\hfill{$\square$}

\medskip
For cube tiling codes $V\subset S^6$ based on Theorem \ref{al} and Corollary \ref{dec} it is possible to prove a weaker version of Theorem \ref{stab5} (we give it without proof):

\begin{tw} 
\label{stab6}
Let $V\subset S^6$ be a cube tiling code. There is a simple component $Q\subset V$ such that it is possible to pass from $V$ to $U(Q)$ by gluing and cutting with $Q$ fixed, where  $U(Q)$ is the simple cube tiling code containing $Q$.
\hfill{$\square$}  
\end{tw}

\subsection{Open problems}

\bigskip
Since every cube tiling code $V\subset S^7$ contains a twin pair (\cite{BHM}), the first problem (it was posed in \cite{DIP}) reads as follows:

\medskip
{\bf Problem 1}. Is it true that every cube tiling codes $V,W\subset S^7$ are strongly equivalent?

\medskip
It is very doubtful that it is possible to resolve this problem along the same lines  as for $d=6$, because the characterization of all pairs $V,W\subset S^6$ strongly equivalent polybox codes $V,W$ such that $|V|\leq 31$ is very difficult. However, it seems that the method applied here allows to reduce the problem to examine polybox codes $V,W\subset S^7$ for $S=\{a,a',b,b',c,c'\}$. Additionally, we may try to establish conditions which has to be satisfied by supposed counterexample (that is, a cube tiling code $V\subset S^7$ which is not strongly equivalent to a fixed simple cube tiling code $U\subset S^7$) and next make suitable computer experiments (such as in \cite{BHM}). For example, by Theorem \ref{gc}, $V$ cannot be strongly equivalent to any cube tiling code $\bar{V}$ such that $\bar{V}=\bar{V}^{i,s}\cup \bar{V}^{i,s'}$ for some $i\in [7]$ and $s\in S$ (compare the proof of Theorem \ref{gc}). 

\medskip
It is interesting question whether the following generalization of Theorem \ref{stab5} holds true:

\medskip
{\bf Problem 2}. Let $V\subset S^d$ be a cube tiling code, and let $Q\subset V$ be a simple component of $V$. Suppose that cube tiling codes $V$ and $U(Q)$ are strongly equivalent, where $U(Q)$ is the simple cube tiling code containing $Q$. Can we pass from $V$ to $U(Q)$ by gluing and cutting with $Q$ fixed? If not, can we characterize cube tiling codes for which such fixed simple components exist?

\medskip
Beside the special pair we know only a few (about 30)  pairs of disjoint and equivalent polybox codes $V,W$ which do not contain twin pairs, and thus which are not strongly equivalent (we do not take into consideration twin pairs free cube tiling codes obtained by the method given in \cite{LS2}). The next problem asks in some sense of the role of twin pairs in equivalent but not strongly equivalent polybox codes:

\medskip
{\bf Problem 3}. Let $V,W\subset S^d$ be equivalent polybox codes which are not strongly equivalent. Suppose that $V$ and $W$ cannot be decomposed such as codes in Corollary \ref{dec}. Can both codes $V$ and $W$ contain twin pairs?

\medskip
To show that two cube tiling codes $V,W\subset S^d$ are strongly equivalent it is enough to show that it is possible to pass from $V$ and $W$ by gluing and cutting to two simple codes $P$ and $Q$, respectively. A path (in the graph $G$ (compare Section 1)) joining $V$ and $W$ and passing through $P$ and $Q$ usually is not the shortest. Moreover, in practice we do not know, beside very small dimensions $d$, the set $\ka T(S^d)$ as it is too large (let us recall that $|\ka T(S^6)|>10^{84}$ for $|S|=64$). Therefore, it is interesting how to find such shortest path (compare \cite{R}):

\medskip
{\bf Problem 4}. 
Let $V,W\in \ka T(S^d)$. Based on the structure of cube tiling codes $V$ and $W$ find a shortest path between $V$ and $W$ in the graph $G$. 



\medskip
The solution of the problem considered in  presented paper has its roots in our investigations on Keller's cube tiling conjecture. There are many other interesting problems dealing the structure of cube tilings and we are believe that they will be examined in the future.

\end{document}